\documentclass{article}

\usepackage{amsmath}
\usepackage{amsfonts}
\usepackage{amscd}
\usepackage{amssymb}

\input xypic

\newtheorem{defn}{Definition}[section]
\newtheorem{thm}[defn]{Theorem}
\newtheorem{prop}[defn]{Proposition}
\newtheorem{lemma}[defn]{Lemma}

\newtheorem{claim}[defn]{Claim}

\newtheorem{rem}[defn]{Remark}

\newtheorem{ass}[defn]{Hypothesis}
\newtheorem{ques}[defn]{Question}

\newcommand{\lm}{\ensuremath{\longrightarrow}}

\newcommand{\eps}{\varepsilon}

\DeclareMathOperator{\Hom}{\mbox{Hom}}

\DeclareMathOperator{\shom}{\ensuremath{\mathcal{H}\mathit{om}}}
\DeclareMathOperator{\sext}{\ensuremath{\mathcal{E}\mathit{xt}}}
\DeclareMathOperator{\sspec}{\ensuremath{\mathcal{S}\mathit{pec}}}
\DeclareMathOperator{\send}{\ensuremath{\mathcal{E}\!\mathit{nd}}}
\DeclareMathOperator{\End}{\mbox{End}}
\DeclareMathOperator{\Aut}{\mbox{Aut}}
\DeclareMathOperator{\Ext}{\mbox{Ext}}

\DeclareMathOperator{\im}{\mbox{im}\,}
\DeclareMathOperator{\spec}{\mbox{Spec}\,}

\DeclareMathOperator{\Mod}{\mbox{Mod}}

\DeclareMathOperator{\Div}{\mbox{Div}\,}
\DeclareMathOperator{\rk}{\mbox{rank}\,}
\DeclareMathOperator{\coker}{\mbox{coker}\,}

\DeclareMathOperator{\Hilb}{\mbox{Hilb}}

\DeclareMathOperator{\qcoh}{\mbox{QCoh}}
\DeclareMathOperator{\mmm}{\mathfrak{m}}
\DeclareMathOperator{\s}{\sigma}

\DeclareMathOperator{\z}{\zeta}
\DeclareMathOperator{\w}{\omega}

\DeclareMathOperator{\PP}{\mathbb{P}}

\DeclareMathOperator{\calm}{\mathcal{M}}
\DeclareMathOperator{\calo}{\mathcal{O}}

\DeclareMathOperator{\bara}{\bar{A}}

\DeclareMathOperator{\ox}{\mathcal{O}_{X}}

\begin{document}

\begin{center}
\LARGE \textbf{Rational curves and ruled orders on surfaces}
\end{center}

\begin{center}
  DANIEL CHAN
\footnote{This project was supported by the Australian Research Council, Discovery Project Grant DP0880143.}
, KENNETH CHAN
\end{center}

\begin{center}
  {\em University of New South Wales}
\end{center}

\begin{center}
e-mail address:{\em danielc@unsw.edu.au,kenneth@unsw.edu.au}
\end{center}

\begin{abstract}
We study ruled orders. These arise naturally in the Mori program for orders on projective surfaces and morally speaking are orders on a ruled surface ramified on a bisection and possibly some fibres. We describe fibres of a ruled order and show they are in some sense rational. We also determine the Hilbert scheme of rational curves and hence the corresponding non-commutative Mori contraction. This gives strong evidence that ruled orders are examples of the non-commutative ruled surfaces introduced by Van den Bergh.
\end{abstract}

Throughout, we work over an algebraically closed base field $k$ of characteristic zero.

\section{Introduction}  \label{sintro}  

Over the last decade, orders on projective surfaces have proved to be a fruitful area of research. As sheaves of algebras on surfaces, they provide examples of non-commutative surfaces which can be studied using the vast arsenal of techniques in algebraic geometry. Notable advances include a non-commutative version of Mori's minimal model program for orders on projective surfaces [CI], as well as an Enriques style classification of them [CI],[CK05]. On a coarse level, the classification gives a trichotomy for the so-called terminal orders without ``exceptional curves'' into: i) del Pezzo orders, ii) (half)-ruled orders and iii) minimal models. Del Pezzo orders have been studied in the generic case since they are examples of Sklyanin algebras (see [ATV90],[ATV91],[A92]), whilst some exotic examples have been studied via the cyclic covering construction (see [CK11] and Lerner's thesis [Ler]). 

In this paper, we wish to study ruled orders which can be described as follows. We start with a smooth projective surface $Z$  ruled over a curve $C$ say by $f:Z \lm C$. Let $A$ be a maximal order of rank $e^2$ on $Z$ such that $k(A) := A \otimes_Z k(Z)$ is a division ring. In keeping with the philosophy of the minimal model program, we shall further assume that $A$ is terminal (see \S2 for a definition), a condition on ramification data which ensures in particular that locally at any closed point, $A$ has global dimension two. If the ramification locus consists of a bisection $D$ together with some fibres, then we will say that $A$ is a {\em ruled order}.

Unlike the complete local structure of orders on surfaces, which has been well-studied (see [A86], [RV], [CI], [CHI]), little is known about the global structure of orders on projective surfaces and we wish in particular to address this lacuna. To guide our study of ruled orders we pose the motivating question, ``In what sense is a ruled order $A$ ruled?''. The approach in [CI] essentially corresponds to the idea that the canonical divisor of a ruled order is positive on fibres. However, it frequently occurs that equivalent concepts in commutative algebraic geometry are inequivalent in non-commutative algebraic geometry. 

Another incarnation of ruledness comes from joint work of the first author with Nyman [CN]. For the purposes of this introduction, a rational curve will be a quotient map $A \lm M$ where $M$ is an $A$-module supported on a smooth rational curve $F$, such that as a sheaf on $F$, i) $M$ is locally free of rank $e$ and ii) $h^0(M) = 1, h^1(M) = 0$. This generalises the notion of line modules on the quantum projective plane (see [ATV91], [A92]). Naturally, we would like to show that ruled orders have lots of rational curves and compute the moduli of them. Furthermore, interest in these rational curves stems from a non-commutative Mori contraction theorem which takes the following form. Given such a rational curve $M$ which is $K$-non-effective and has self-intersection zero (see \S2 for definitions), there is a morphism of non-commutative schemes $\sspec_Z A \lm Y$ where $Y$ is the Hilbert scheme of deformations of $A \lm M$ (see theorem~\ref{tCN} for an explanation of  terms and details). Furthermore, $Y$ is a generically smooth (commutative) curve. A natural question thus is
\begin{ques} \label{qMori}  
Does a ruled order have a $K$-non-effective rational curve with self-intersection zero? If so, what are the corresponding component of the Hilbert scheme and non-commutative Mori contraction?
\end{ques}
A natural place to look for rational curves is to consider fibres $F$ of the fibration $f:Z \lm C$ and to consider rank $e$ locally free quotients of $\bara:= A\otimes_Z \calo_F$ with the desired cohomology. This suggests the next
\begin{ques}  \label{qfibre}  
Describe the fibres $\bara$ of $A$. In particular, is $\bara$ ``rational'' in the sense that $H^1(\bara) = 0$? 
\end{ques}
Now if $F$ is not a ramification curve, then $\bara$ is an order and so embeds in a maximal order, which by Tsen's theorem has the form $\send_F V$ for some vector bundle on $V$ on $F$. Thus one need only determine possibilities for $V$ and subalgebras of $\send_F V$. Question~\ref{qfibre} is nevertheless fairly subtle. Indeed, ramification data, and hence the notion of a ruled order is Morita invariant whilst cohomology is not. To select the best ``model'' in the Morita equivalence class, we will assert the minimality of the second Chern class $c_2(A)$. A key tool is the notion of Morita transforms, already introduced by Artin-de Jong in the case of Azumaya algebras, to construct Morita equivalences. 

The answer to question~\ref{qfibre} is nicest when $F$ is a fibre which meets the ramification divisor transversally (necessarily then in two distinct points). The following comes from theorem~\ref{theredfibres}.
\begin{thm}  \label{tintrohered}  
Let $A$ be a ruled order with minimal $c_2$ in its Morita equivalence class. Let $F \subset Z$ be a fibre which meets the ramification locus transversally in two distinct points. Then 
$$\bara:= A \otimes_Z \calo_F \simeq 
\begin{pmatrix}
\calo & \calo(-1) & \cdots & \calo(-1) \\
\calo(-1) & \calo & \ddots & \vdots \\
\vdots & \ddots & \ddots & \calo(-1) \\
\calo(-1) & \cdots & \calo(-1) & \calo
\end{pmatrix}.
$$
In particular, $H^1(\bara) = 0$.
\end{thm}

This theorem shows there are $e$ rational curves which arise as quotients of $\bara$, namely, the columns of the matrix form above. The natural way to enumerate them is as the points of $\spec H^0(\bara)$. That suggests that one can get a handle on the Hilbert scheme of rational curves by considering the following setup. Note that by proposition~\ref{pf*A}, $f_* A$ is a commutative coherent sheaf of algebras on $C$. There is also a natural multiplication map $\Psi: A \otimes_Z f^*f_*A \lm A$ which respects the left $A$-module structure and right $f^*f_*A$-modules structure. 

\begin{thm}  \label{tmain}  
Let $A$ be a ruled order with minimal $c_2$ in its Morita equivalence class. Suppose further that every ramified fibre has ramification index equal to $\deg A$ (e.g. when $\deg A$ is prime). Then 
\begin{enumerate}
\item $C':=\sspec_C f_*A$ is a smooth projective curve. 
\item $C'$ is a component of the Hilbert scheme $\Hilb A$ and $\Psi: A \otimes_Z f^*f_*A \lm A$ is the universal rational curve.
\item The non-commutative Mori contraction is given by the ``morphism of algebras'' $f_*A \lm A$.  
\end{enumerate}
\end{thm}

This is the main theorem and the proof involves analysing the fibres $\bara$ of $A$ case by case. The paper has been organised as follows. Section~\ref{ssetup} contains background material. We deal with the generic fibre in section~\ref{sgeneric}. To obtain nice results at other fibres, we will have to impose minimality of the second Chern class and so in section~\ref{smorita}, we study the procedure of Morita transforms for constructing Morita equivalences. In particular, we compute formulas for how $c_2$ changes with Morita transforms. In section~\ref{shered}, we study hereditary fibres, namely, those where $F$ intersects the ramification locus transversally in 2 distinct points. In particular, we prove theorem~\ref{tintrohered}. To study the other fibres, we embed them in a minimal hereditary order $\bara_1$ and in section~\ref{sheredover}, we give a description of $\bara_1$ under some hypotheses. Sections~\ref{sfibretan},\ref{sfibrenode},\ref{sramfibres} deal with the three possible non-hereditary fibres, namely, i) $F$ is tangential to the bisection $D$, ii) $F$ contains a node of $D$ and iii) $F$ is a ramification curve. The key tool to analysing these cases is the notion of rational filtrations introduced in section~\ref{sratgood}. The construction of the requisite rational filtration involves the over-order $\bara_1$. Our analysis shows that with the hypotheses of the theorem, we do indeed have $H^1(\bara) = 0$. Finally in section~\ref{smori}, we give the proof of theorem~\ref{tmain}. 

Presumably the results of theorem~\ref{tmain} hold for all ruled orders and extending results to the general setting is the subject of ongoing work. The problem with ramified fibres is that in this case, $\bara$ is not an order so our methods do not apply.

\section{Background and Setup}  \label{ssetup}  

The theory of the birational classification of orders on projective surfaces developed in [CI] mirrors that of commutative surfaces. We recall it here in this section as well as the theory of non-commutative Mori contractions in [CN]. 

Let $Z$ be a smooth projective surface. An {\em order} on $Z$ is a torsion-free coherent sheaf of algebras $A$ on $Z$ such that $k(A):= A \otimes_Z k(Z)$ is a central simple $k(Z)$-algebra. We can thus consider orders as subalgebras of a central simple algebra, and so order them by inclusion. The maximal orders (with respect to inclusion) correspond to the normal surfaces. Orders are Azumaya on some dense open subset $U \subseteq Z$. Its complement $Z-U$ is called the {\em ramification locus}. It is a union of {\em ramification curves} (see [CI, section~2.2] for more details on ramification).

The role of smooth models in the birational classification of orders on surfaces is played by terminal orders, which can be defined as follows. We say a maximal order $A$ on $Z$ is {\em terminal} if 
\begin{enumerate}
\item the ramification locus is a normal crossing divisor, and 
\item if $p$ is a node of the ramification locus lying in the intersection of ramification curves $C_1,C_2$, then the ramification of $A$ at $C_i$ is itself totally ramified at $p$ for at least one of the curves $C_1$ or $C_2$. 
\end{enumerate}
Further details on condition~ii) of the definition can be found in [CI, section~2.2]. The complete local structure of terminal orders has been determined in [CI, section~2.3]. To describe this we use the skew-power series ring $k_{\zeta}[[x,y]]:= k\langle\langle x,y \rangle\rangle/(yx - \zeta xy)$ where $\zeta \in k$ is a root of unity. 
\begin{thm} \label{tlocalterm}  
Let $A$ be a terminal order on $Z$ and $p \in Z$ be a closed point. Then identifying $\hat{\calo}_{Z,p}$ with $k[[u,v]]$ appropriately, we have a $k[[u,v]]$-algebra isomorphism of $A \otimes_Z \hat{\calo}_{Z,p}$ with a full matrix algebra in 
$$ \begin{pmatrix}
 k_{\zeta}[[x,y]]& \cdots & \cdots & k_{\zeta}[[x,y]] \\
(x) & k_{\zeta}[[x,y]] & \ddots & \vdots \\
\vdots & \ddots & \ddots & \vdots \\
(x) & \cdots & (x) & k_{\zeta}[[x,y]]
\end{pmatrix}
\subseteq M_n(k_{\zeta}[[x,y]])$$
for some primitive $e$-th root of unity $\zeta$. Here $x^e = u, y^e = v$ and the ramification locus is $uv=0$. The ramification indices along $u=0$ and $v=0$ are $ne$ and $e$ respectively. 
\end{thm}
If one wants, one can use this description as a definition for terminal orders. That is, a maximal order on a smooth projective surface is terminal if and only if complete locally at any closed point, it has the form in the theorem. 

To simplify the treatment from now on, we will almost exclusively consider orders in a division ring only, that is, those such that $k(A)$ is a division ring. We will justify this assumption shortly.

One of the main results in [CI], is a theory of Mori contractions for terminal orders on projective surfaces. These Mori contractions were classified and described in complete analogy with the commutative case. In particular, the Mori contraction corresponding to a ruled surface has non-commutative counterparts dubbed ruled orders. To define them, we start with a ruled surface $f:Z \lm C$ and let $A$ be a terminal order on $Z$ in a $k(Z)$-central division ring $k(A)$. We say that $A$ is a {\em ruled order} if furthermore, the ramification of $A$ is the union of a bisection $D$ of $f$ with some fibres. In this case, the {\em degree} say $e:= \sqrt{\rk_Z A}$ of $A$ is the same as the ramification index of $A$ on $D$ by the Artin-Mumford sequence (see [AM, theorem~1]). 

Following [AZ94], we consider a non-commutative scheme to be a $k$-linear abelian category $\mathcal{C}$ together with a distinguished object $\calo$ called the {\em structure sheaf}. Given an order $A$ we can obtain such a non-commutative scheme by first setting $\mathcal{C}=A-\Mod$ the category of quasi-coherent $A$-modules. For the choice of the structure sheaf, we usually pick $\calo = A$ itself. This means that the cohomology of $A$-modules will be their usual Zariski cohomology as sheaves on $Z$. We will allow ourselves to change $\calo$ to some local pro-generator $P$ which has the same effect as replacing our order $A$ with the Morita equivalent one $A':= \send_Z P$ and ``restoring'' the structure sheaf to $\calo = A'$.
Now any terminal order in a central simple $k(Z)$-algebra is Morita equivalent to a terminal order in a division ring, which explains why we restricted our attention to this special case. We will have further need of Morita equivalences later. We will sometimes write $\sspec_Z A$ for the non-commutative scheme $(\mathcal{C}, \calo)=(A-\Mod, A)$. 

The notion of non-commutative Mori contractions in [CN] is very different from the Mori contractions for terminal orders in [CI]. To describe the former, we need to introduce some terminology. As motivation, recall that a ruled surface corresponds to the Mori contraction of a $K$-negative rational curve with self-intersection zero. Let $A$ be a terminal order on a smooth projective surface $Z$. A {\em rational curve} on $\sspec_Z A$ is a cyclic (left) $A$-module $M$ which is pure of dimension one and has cohomology $H^0(M) = k, H^1(M) = 0$. We say that the rational curve $M$ is {\em $K$-non-effective} if $H^0(\omega_A \otimes_A M) = 0$ where $\omega_A$ is the $A$-bimodule $\shom_Z(A, \omega_Z)$ (see [CK03, section~3, proposition~5] for more information about $\omega_A$). Finally we say that the rational curve $M$ has self-intersection zero if the Euler form 
$$ \chi(M,M) = \sum_{i=0}^2 \dim_k \Ext_A^i(M,M)=: -M^2 = 0.$$ 

Recall from [AZ01, section~E], that the Hilbert scheme $\Hilb A$ parametrising quotients of $A$ is a separated locally finite type scheme which is a countable union of projective schemes. We may thus consider the component $Y$ of the Hilbert scheme containing the point corresponding to a rational curve $M$ and let $\calm$ be the universal family on $Y$. We call $Y$ the {\em Hilbert scheme of deformations of $M$}. We summarise some results from [CN] in the special case of terminal orders on surfaces. Below, a {\em 1-critical} module is a 1-dimensional module all of whose non-trivial quotients have dimension $< 1$. 

\begin{thm} \label{tCN}  
Let $A$ be a terminal order on a smooth projective surface $Z$ and $M$ a 1-critical $K$-non-effective rational curve with self-intersection zero. Let $Y$ be the Hilbert scheme of deformations of $M$ and $\calm$ be the universal family. Then 
\begin{enumerate}
\item $Y$ is a projective curve which is smooth at the point corresponding to $M$.
\item If $\pi_Y: Y \times Z \lm Y, \pi_Z: Y \times Z \lm Z$ denote projection maps then the Fourier-Mukai type transform
$$  \pi_{Z*}(\calm \otimes_{Y \times Z} \pi^*_Y(-)): \qcoh(Y) \lm A-\Mod  $$
is an exact functor with a right adjoint.
\end{enumerate}
\end{thm}
\textbf{Remark:} A morphism of non-commutative schemes is just a pair of adjoint functors so we will refer to the Fourier-Mukai transform in ii) as the {\em non-commutative Mori contraction contracting $M$}.

\noindent
\textbf{Proof.} One checks easily that $\sspec_Z A$ is a non-commutative smooth proper surface in the sense of [CN, \S~3] (the Serre functor is given by $\omega_A \otimes_A-$ and the dimension function is the usual dimension of coherent sheaves). We may thus invoke [CN, corollary~9.7] to obtain i) and [CN, theorem~10.2] to obtain ii).
\vspace{2mm}

We would like to show that ruled orders give examples of non-commutative Mori contractions and examine closely the nature of these contractions. The following elementary result is crucial.
\begin{prop} \label{pf*A}  
Let $A$ be a ruled order on the ruled surface $f:Z \lm C$. Then the $\calo_C$-algebra $f_*A$ is commutative and in fact, is a sheaf of integral domains on $C$.
\end{prop}
\textbf{Proof.} First note that $f_*A$ is a torsion-free coherent sheaf on $C$ so the subring $f_*A \otimes_C k(C) \subset k(A)$ is a division ring. Tsen's theorem shows that in fact $f_*A \otimes_C k(C)$ is a field so in particular $f_*A$ is commutative.
\vspace{2mm}

We may thus consider the finite cover $C':= \sspec_C f_*A$ of $C$. Note that $Z':= C' \times_C Z = \sspec_Z f^*f_* A$ and that $A$ is an $(A, f^*f_*A)$-bimodule and hence, also can be considered as a family of $A$-modules over $C'$. Furthermore, there is a natural map
$$ \Psi: A \otimes_Z f^*f_*A \lm A $$ 
induced by multiplication which is surjective since $f^*f_*A \supset \calo_Z$. Hence if $C'_{\text{fl}} \subseteq C'$ is the locus where $A$ is flat, then there is an induced map $h:C'_{\text{fl}} \lm \Hilb A$. We wish of course to show that $C'_{\text{fl}} = C'$ and $h$ is an isomorphism of $C'$ with a Hilbert scheme of deformations of a rational curve. The former would follow easily if we knew {\em a priori} that $C'$ were smooth, but unfortunately, we will only be able to prove smoothness by first showing the Hilbert scheme is smooth and $h$ is well-defined on all of $C'$. 

\section{Generic behaviour}  \label{sgeneric}  

Consider a degree $e$ ruled order $A$ and the map $\Psi:A \otimes_Z f^*f_*A \lm A$ introduced in section~\ref{ssetup}. In this section, we study the behaviour of $\Psi$ generically on $C$. We start by studying $f_*A$ and more generally $R^if_* A$. First note that $A$ is locally free as a sheaf on $Z$ so is flat over $C$. The following result has already been noted by Artin and de Jong in [AdJ, proposition~4.3.1].

\begin{prop} \label{pfAgen}  
We have $\rk f_*A = e$ and $R^1f_* A$ is a torsion sheaf.
\end{prop}
\textbf{Proof.} Let $\eta \in C$ be the generic point and $A_{\eta}$ be the pullback of $A$ to $Z_{\eta}:=Z \times_C \eta$. Since cohomology commutes with the flat base change $\eta \hookrightarrow C$, it suffices to show that $\dim_{k(C)} H^0(A_{\eta}) = e, H^1(A_{\eta})=0$.

We first use a formula of Artin-de Jong for the Euler characteristic (see [AdJ, (4.1.2)] or [Chan]) which is easily derived from examining the \'etale local forms of an order.
$$ \frac{\chi(A_{\eta})}{(\deg A_{\eta})^2} = \chi(\calo_{Z_{\eta}}) - \frac{1}{2}\sum_{i=1}^s (1 - \frac{1}{e_i})$$
where $e_1,\ldots, e_s$ are the ramification indices of $A_{\eta}$ written with multiplicity. Here the multiplicity of the ramification index $e_i$ corresponding to the ramification point $p_i \in Z_{\eta}$ is $[k(p_i):k(C)]$. In our case, $A$ is ramified on a bisection union some fibres, so either $A_{\eta}$ is ramified on two distinct points with multiplicity one, or one point with multiplicity two. In either case, the ramification indices with multiplicity are $e_1=e_2 = e$. Also, $Z$ is ruled so $\chi(\calo_{Z_{\eta}}) = 1$. Thus $\chi(A_{\eta}) = e$. Note that $H^0(A_{\eta})(z) \simeq H^0(A_{\eta})k(Z) \subset k(A)$ is a commutative subalgebra of $k(A)$ so $\dim_{k(C)} H^0(A_{\eta}) \leq e$. The only possibility is $\dim_{k(C)} H^0(A_{\eta}) = e, H^1(A_{\eta})=0$.
\vspace{2mm}

It will become apparent later, that $R^1f_*A$ is not necessarily zero. The vanishing locus of $R^1f_*A$ will be the locus where the map $\Psi$ is well-behaved. Let $c \in C$ be a closed point and $A_c$ denote the restriction of $A$ to the fibre $f^{-1}(c) \subset Z$. By cohomology and base change results [Hart, theorem~III.12.11], we have $R^1f_*A \otimes_C k(c) \simeq H^1(A_c)$ and if this is zero, we also have $f_*A \otimes_C k(c) \simeq H^0(A_c)$ which we note is commutative. 

Suppose that the fibre $F:=f^{-1}(c)$ is not a ramification curve. Then $A_c$ is an order on $F \simeq \PP^1$. If furthermore the fibre intersects the ramification locus transversally, then by the local structure theory of terminal orders given in theorem~\ref{tlocalterm}, $A_c$ is an hereditary order ramified at two points with ramification indices both $e$. 
The splitting theorem of [C05], ensures that the left $A_c$-module $A_c$ is the direct sum of $e$ locally projective $A_c$-modules, say $P_1,\ldots , P_e$, each of rank $e$. The image of the identity $1\in A_c$ in each summand $P_i$ is an idempotent. This gives the useful 
\begin{lemma}  \label{lheredPeirce} 
An hereditary order $\bara$ on $\PP^1$ which is ramified at two or fewer points has a Peirce decomposition 
$$\bara = \bigoplus_{i,j = 1}^{e} \calo_{\PP^1}(p_{ij}) $$
for some divisors $p_{ij} \in \Div \PP^1$. 
\end{lemma}

We may now describe the behaviour of $\Psi:A \otimes_Z f^*f_*A \lm A$ for generic values of $c$.

\begin{thm}  \label{tgeneric}  
Let $c \in C$ be a closed point such that $f^{-1}(c)$ intersects the ramification locus of $A$ transversally and furthermore, $H^1(A_c) = 0$. Also, let $C^0 \subseteq C$ be the dense open set of such points. Then
\begin{enumerate}
\item We have
$$A_c \simeq 
\begin{pmatrix}
\calo & \calo(-1) & \cdots & \calo(-1) \\
\calo(-1) & \calo & \ddots & \vdots \\
\vdots & \ddots & \ddots & \calo(-1) \\
\calo(-1) & \cdots & \calo(-1) & \calo
\end{pmatrix}.
$$
\item $\pi:C' = \sspec_C f_*A \lm C$ is \'etale over $C^0$.
\item Suppose $c' \in C'$ is a point in the pre-image of $c \in C^0$. Then $A\otimes_{C'} k(c')\simeq \calo_{\PP^1} \oplus \calo_{\PP^1}(-1)^{\oplus (e-1)}$. In particular, $A$ is flat over $\pi^{-1}(C^0)$. 
\item For $c' \in \pi^{-1}(C^0)$, the $A$-module $A\otimes_{C'} k(c')$ is a 1-critical $K$-non-effective rational curve with self-intersection zero.
\end{enumerate}
\end{thm}
\textbf{Proof.} To see part i), note that $H^1(A_c) = 0$ implies $\deg p_{ij} \geq -1$ for all $i,j$. Also, $p_{ii} = 0$ so $H^0(A_c) = e$ forces in fact $\deg p_{ij} = -1$ for all $i \neq j$. Let $\eps_1,\ldots,\eps_e$ be the $e$ diagonal idempotents in i). Then $f_*A \otimes_C k(c) = H^0(A_c) = k\eps_1 \times \ldots \times k\eps_e$. Thus the pre-image of $c$ in $C'$ consists of $e$ distinct points so $\pi$ is \'etale over $C^0$. 

To prove iii), note that $c'$ corresponds to an idempotent of $H^0(A_c)$, say $\eps_i$. Then $\Psi \otimes_{C'} k(c')$ is just the multiplication map 
$$ A_c \otimes_k k\eps_i \lm A_c \eps_i .$$
Now $A_c\eps_i$ is just the $i$-th column of $A_c$ so is isomorphic to $\calo_{\PP^1} \oplus \calo_{\PP^1}(-1)^{\oplus (e-1)}$. Flatness now follows from the fact that the Hilbert polynomial is constant (see [Pot, theorem~4.3.1] or [Hart, theorem~III.9.9]).

We prove iv). Let $M = A \otimes_{C'} k(c')$ which we note is certainly 1-critical. Also $h^0(M) = 1$ and $h^1(M) = 0$ so $M$ is indeed a rational curve. Now if $D$ denotes the bisection that $A$ is ramified on, then [CK, \S~3, proposition~5] shows that we can write the canonical bimodule $\w_A$ in the form $I \otimes_Z \calo(F)$ where $F$ is a linear combination of fibres and $I$ is an ideal of $A$ such that $I^e = A \otimes_Z \calo(-D-F')$ for some $F'\geq 0$ a linear combination of fibres. It follows that $\w_A \otimes_A M \simeq IM$. Note that $M$ is generated as an $A$-module by the unique copy of $\calo$ in $M$, so if $h^0(IM) \neq 0$ then $IM = M$, a contradiction since $I^eM \simeq M(-D)$. This shows that $M$ is $K$-non-effective. We need to show that the Euler form $\chi(M,M) = 0$. Now $\Hom_A(M,M) = k$ and Serre duality gives $\Ext^2_A(M,M) = 0$ so we need only show that $\Ext^1_A(M,M) = k$. Note first that for $c_1 \in C$ distinct from $c$ we have $\chi(M,A_c) = \chi(M,A_{c_1}) = 0$ since $M$ and $A_{c_1}$ have disjoint support. Serre duality shows that $\Ext^2_A(M,A_c) = 0$ while direct computation from the matrix form in i) shows that $\Hom_A(M,A_c) = k$. Hence $\Ext^1_A(M,M)$ is a direct summand of $\Ext^1_A(M,A_c) = k$. We are thus reduced to showing that $\dim_k \Ext^1_A(M,M) \geq 1$. This follows since the family $A/C'$ gives non-trivial flat deformations of $M$. 
This completes the proof of the theorem.
\vspace{2mm}

The above proof can now be used to cover some cases where $f^{-1}(c)$ does not intersect the ramification locus transversally.

\begin{prop}  \label{pflatcrit}  
Suppose $c \in C$ is a closed point with $H^1(A_c) = 0$ and $c' \in C'$ lies over $c$. Let $M = A \otimes_{C'} k(c') = A_c \otimes_{H^0(A_c)} k(c')$. If $M$ is locally free of rank $e$ and $\chi(M) = 1$ then $A/C'$ is flat at $c'$ and $M$ is a $K$-non-effective rational curve with self-intersection zero. 
\end{prop}
\textbf{Proof.} 
We know that $M$ is a quotient of $A_c$ so $H^1(A_c) = 0$ implies that $H^1(M) = 0$. Thus $H^0(M) = k$ so $M$ is indeed a rational curve. In fact, $M\simeq \calo_{\PP^1} \oplus \calo_{\PP^1}(-1)^{\oplus (e-1)}$ so as before, we see that $A/C'$ is flat over $c'$. Continuity of Euler characteristic shows that $M$ has self-intersection zero. The proof of $K$-non-effectivity in the theorem shows that $M$ is also $K$-non-effective. 

It will be convenient to introduce the following
\begin{defn}  \label{dgood}  
We say that the closed fibre $F = f^{-1}(c) \subset Z$ is {\em good} if $H^1(A_c) = 0$ and the conclusions of theorem~\ref{tgeneric}iii),iv) hold. By proposition~\ref{pflatcrit}, $F$ is good if $H^1(A_c) = 0$ and for each $c' \in \pi^{-1}(c)\subset C'$ we have that $A\otimes_{C'} k(c')$ is a 1-critical $A$-module whose underlying sheaf is isomorphic to $\calo_{\PP^1} \oplus \calo_{\PP^1}(-1)^{\oplus (e-1)}$. The {\em good locus} is the set $C_{\text{good}} \subseteq C$ of such points $c$ where $f^{-1}(c)$ is good.
\end{defn}
Recall that $C'_{\text{fl}}$ denotes the locus where $A$ is flat over $C'$. The next result shows how good fibres are indeed good for our purposes.
\begin{prop}  \label{pgood}  
The dense open set $C^0\subseteq C$ of theorem~\ref{tgeneric} is contained in $C_{\text{good}}$ so in particular $C_{\text{good}}$ is dense open. Also, $\pi^{-1}(C_{\text{good}}) \subseteq C'_{\text{fl}}$. The map $\Psi:A \otimes_Z f^*f_* A \lm A$ induces an open immersion $h:\pi^{-1}(C_{\text{good}}) \lm \Hilb A$. Furthermore, $\pi^{-1}(C_{\text{good}})$ is smooth. 
\end{prop}
\textbf{Proof.} The theorem shows $C^0 \subseteq C_{\text{good}}$. Note that $A$ is flat over $\pi^{-1}(C_{\text{good}})$ since the Hilbert polynomial is constant so $\pi^{-1}(C_{\text{good}}) \subseteq C'_{\text{fl}}$ and the map $h$ is well-defined. Furthermore, theorem~\ref{tCN} ensures that $\Hilb A$ is smooth 1-dimensional at all points of $\im h$. It thus suffices now to show that $h$ is generically injective. The easiest way to see this is to consider a point $c \in C^0$ and the matrix form for $A_c$ given in theorem~\ref{tgeneric}. Then the points of $\pi^{-1}(c)$ correspond to the columns of this matrix form which correspond to distinct points of $\Hilb A$. 
\vspace{2mm}

Naturally, we would like to show that all fibres are good.

\section{Morita transforms}  \label{smorita}  

The previous section shows that the condition $H^1(A_c) = 0$ is extremely important to make proofs work. To attain this condition, we pass to a Morita equivalent order. The question of which order is a little delicate. One way to restrict the choices is to insist that $c_2(A)$ is minimal. We note that by results of Artin-de Jong [AdJ], there is a lower bound for $c_2$ for maximal orders in a fixed division ring. In the case of terminal orders, any two such orders are Morita equivalent by [CI, corollary~2.13].

We describe here the method of Morita transforms for changing to a Morita equivalent order in the same central simple algebra. This was called ``elementary transformations'' by Artin and de Jong, a term we have avoided since it may have an alternate meaning when dealing with ruled orders. Artin-de Jong's treatment was restricted to Azumaya algebras (see [AdJ, \S~8.2] or [dJ, section~2]).

We thank Michael Artin and Johan de Jong for allowing us to include the following result and their proof.
\begin{prop}  \label{pboundc2}  
Let $Z$ be a smooth projective surface and $Q$ be a $k(Z)$-central division ring of degree $e$. Given these data, there is a constant $\gamma$ such that for any maximal order $A$ in $Q$ we have $\chi(A) \leq \gamma$. In particular, there is a lower bound on the second Chern class $c_2(A)$ depending only on $Z$ and $Q$.
\end{prop}
\textbf{Proof.} We quote almost verbatim from [AdJ, \S~7.1]. From the \'etale local form of maximal orders at codimension one points, $c_1(A)$ is independent of the choice of maximal order (see [AdJ, \S~7.1] or [Chan]). 
Hence we need only bound the Euler characteristic. Consider the reduced determinant map $\det:A \lm \calo_Z$. If $L$ is a line bundle then this map extends to a morphism $\det:A \otimes_Z L \lm L^{\otimes e}$, and hence gives a map of affine varieties
$$H^0(A \otimes_Z L) \lm H^0(L^{\otimes e}) .$$
Now $\det$ is multiplicative and $Q$ a division ring, so the fibre $\det^{-1}(0)$ consists of exactly one point, namely $0$. Chevalley's theorem then shows
\begin{equation} 
\dim_k H^0(A \otimes_Z L) \leq \dim_k H^0(L^{\otimes e}) \label{ec2} 
\end{equation}
In particular, $h^0(A) \leq h^0(\calo_Z)$. We now bound $h^2(A)$. By Serre duality we have $h^2(A) = h^0(A^* \otimes_Z \omega_Z)$. We may embed $A^*$ in $D$ (for example using the reduced trace map) and so assume that $A^* \subset A \otimes_Z \calo_Z(\Delta)$ for a sufficiently ample divisor $\Delta$. Then equation~(\ref{ec2}) shows that 
$$h^2(A) \leq h^0(\omega_Z^{\otimes e}(e \Delta)).$$
The result follows.

\noindent
\textbf{Remark:} From a strictly logical point of view, we will only need the result concerning the upper bound on $\chi(A)$. The result on $c_2(A)$ has been included since it provides the best context for presenting results.
\vspace{2mm}

Consider now a terminal order $A$ on a smooth surface $Z$ and $F \subset Z$ a smooth curve. Let $\bara := A|_F$ and $I< A$ be a left ideal containing $A(-F)$ such that $M:= A/I$ is an $\bara$-module which is pure of dimension one in the sense that no submodules are supported at closed points. Thus $I$ is reflexive. In fact more can be said.

\begin{prop}  \label{pprogen}  
Let $A$ be a degree $d$ terminal order on a quasi-projective surface $Z$. Let $P$ be an $A$-module which is reflexive as a sheaf on $Z$ and with the same rank as $A$. Then $P$ is a local progenerator which induces a Morita equivalence between $A$ and the maximal order $A' := \send_A P$. In fact, for any closed $p \in Z$ we have an isomorphism of completions $\hat{P}_p \simeq \hat{A}_p$.
\end{prop}
\textbf{Proof.} First note that $Z$ is smooth so $P$ is locally free as a sheaf on $Z$. Hence $P$ is locally projective by [Ram, proposition~3.5]. Note that $k(Z) \otimes_Z P \simeq k(A)$ so that $A'$ is also an order in $k(A)$. 

Now $P$ is a reflexive sheaf so the same is true of $A'$. By Auslander-Goldman's criterion [AG, theorem~1.5], maximality of $A'$ will follow from maximality at all codimension one points. Let $C$ be a prime divisor. We check maximality at $C$ by computing the discriminant ideal $\mmm_C^j$ there, where $\mmm_C\triangleleft \calo_{Z,C}$ is the maximal ideal of the local ring at the generic point of $C$. Let $e$ be the ramification index of $A$ at $C$, $f=d/e$ and $p \in C$ be a general point so that either $A$ is Azumaya at $p$ or $p$ lies on the smooth locus of the ramification. Then by theorem~\ref{tlocalterm} we know 
$\hat{A}_p \simeq S^{f \times f}$ where 
$$
S \simeq 
\begin{pmatrix}
k[[x,y]] & \cdots & \cdots & k[[x,y]] \\
(x)      & k[[x,y]] &       & \vdots \\
\vdots   & \ddots & \ddots  & \vdots \\
(x)      & \cdots & (x)     & k[[x,y]]
\end{pmatrix}.$$
Here, we may identify $k[[x,y]]$ with the complete local ring $\hat{\calo}_{Z,p}$ and $x=0$ is a local equation for $C$. Note that there are $e$ distinct indecomposable projective $\hat{A}_p$-modules, say $P_1, \ldots, P_e$ corresponding to the columns of $S$ above, and that 
$$\hat{A}_p \simeq P_1^f \oplus \ldots \oplus P_e^f .$$
We seek to show that $\hat{P}_p \simeq \hat{A}_p$ from which it will follow that the discriminant ideals of $A$ and $A'$ are equal and hence, $A'$ will also be maximal. Now we can write 
$$\hat{P}_p \simeq P_1^{f_1} \oplus \ldots \oplus P_e^{f_e} $$
where $\sum f_i = d$. The discriminant ideal of $A'$ is $\mmm_C^j$ where
$$ j = \frac{1}{2} \Bigl[ d(d - 1) - \sum_{i=1}^e f_i(f_i-1)\Bigr] = 
\frac{1}{2} \bigl( d^2 - \sum_{i=1}^e f_i^2 \bigr) .$$
Now $j$ is maximised precisely when all the $f_i$ equal $f$. On the other hand, $\mmm_C^j$ is contained in the discriminant ideal of the maximal order containing $A'$. Hence all the $f_i$ are equal and $\hat{P}_p \simeq \hat{A}_p$. 

To finish the proof of the proposition, we need only show that $\hat{P}_p \simeq \hat{A}_p$ even for points $p$ at the nodes of the ramification locus. This follows using the same computation as above, but using the complete local structure of terminal orders at such nodes. 

\vspace{2mm}
We call $A':= \send_A I$ the {\em Morita transform of $A$ associated to $M$}. We wish to examine how $c_2$ changes with Morita transforms. Following Artin-de Jong, we consider the subalgebra $B = A \cap A'$ which we write as $\send_A (I \subset A)$, the sheaf of endomorphisms of $A$ which stabilise $I$. From this point of view, we obtain the exact sequence
$$ 0 \lm B \lm A \xrightarrow{\phi} \shom_A(I,M) .$$
Now right multiplication by any section of $A$ clearly sends $A(-F)$ into $A(-F)$ so $\phi$ factors through $A \lm \shom_{\bara}(M',M)$ where $M':=I/A(-F)$. We have thus an exact sequence 
$$ 0 \lm B \lm A \lm \shom_{\bara}(M',M) \lm T \lm 0 $$
for some sheaf $T$ supported on $F$.

We can also write $B= \send_A(A(-F) \subset I)$ from which we see there is also an exact sequence of the form 
$$ 0 \lm B \lm A' \lm \shom_{\bara}(M(-F),M') \lm T' \lm 0 $$
for some sheaf $T'$ supported on $F$. 

\begin{prop} \label{pMorita}  
With the above notation we have
\begin{enumerate}
\item $T= \sext^1_{\bara}(M,M),\ T'=\sext^1_{\bara}(M',M')$.
\item $c_2(A')-c_2(A) = \chi(\shom_{\bara}(M',M))) - \chi(\shom_{\bara}(M(-F),M')) -  \chi(T) + \chi(T')$. 
\end{enumerate}
\end{prop}
\textbf{Proof.} We calculate $T'$ as follows. Let $\bar{I} = I \otimes_Z \calo_F$ which we note is a locally projective $\bara$-module. We have the exact sequence 
$$ 0 \lm M(-F) \lm \bar{I} \lm M' \lm 0 .$$
Now the map $A' \lm \shom_{\bara}(M(-F),M')$ factors through the natural  quotient map
$$A' \lm A' \otimes_Z \calo_F = \send_{\bara}\bar{I}.$$
The formula for $T'$ follows now from the exact sequences
$$ \send_{\bara} \bar{I} \lm \shom_{\bara}(\bar{I},M') \lm \sext^1_{\bara}(\bar{I},M(-F)) = 0 $$
$$ \shom_{\bara}(\bar{I},M')  \lm \shom_{\bara}(M(-F),M') \lm \sext^1_{\bara}(M',M') \lm \sext^1_{\bara}(\bar{I},M') = 0 .$$
The computation for $T$ is easier and uses the same method so will be omitted. 

Part ii) follows from the fact that the rank and first Chern classes of $A$ and $A'$ are the same so Riemann-Roch (see [Pot, \S91, p.154]) gives $c_2(A')-c_2(A) = \chi(A) - \chi(A')$.

\vspace{2mm}

\begin{prop}  \label{pnoExt}  
Suppose further that $F$ is not a ramification curve so that $\bara$ is an order. Then 
\begin{enumerate}
\item With the above notation, $\chi(T) = \chi(T')$.
\item If $A$ has minimal $c_2$ in its Morita equivalence class, then 
$$ \chi(\shom_{\bara}(M',M))) \geq \chi(\shom_{\bara}(M(-F),M')) .$$
\end{enumerate}
\end{prop}
\textbf{Proof.} In view of proposition~\ref{pMorita}, it suffices to prove part~i). Now $M,M'$ are torsion-free sheaves on $F$ so are locally projective $\bara$-modules everywhere except possibly where the ramification locus of $A$ intersects $F$ with multiplicity $\geq 2$. We may thus calculate the Ext-sheaves $T,T'$ by going to the complete local rings at such points. We apply two long exact sequences in cohomology related to the exact sequence
$$ 0 \lm M' \lm \bara \lm M \lm 0 .$$
One application shows that $\sext^1_{\bara}(M',M') \simeq \sext^2_{\bara}(M,M')$. Another application gives the exact sequence
$$ \sext^1_{\bara}(M,\bara)  \lm \sext^1_{\bara}(M,M) 
\lm \sext^2_{\bara}(M,M')  \lm \sext^2_{\bara}(M,\bara) .$$
Now for any $p \in F$, the completion $A \otimes_Z \hat{\calo}_{Z,p}$ is regular of dimension two [CI, theorem~2.12] so $\bara_p := \bara \otimes_F \hat{\calo}_{F,p}$ has injective dimension one. Thus $\sext^2_{\bara}(M,\bara) = 0$. Also, $M$ is torsion-free and hence, a first syzygy so also $\sext^1_{\bara}(M,\bara)=0$. This completes the proof of the proposition. 

\begin{prop} \label{pminc2ruled}  
Let $A$ be a degree $e$ ruled order with minimal $c_2$ in its Morita equivalence class.  Let $F \subset Z$ be a fibre of the ruling which is not a ramification curve and consider an exact sequence of the form
$$ 0 \lm M' \lm \bara \lm M \lm 0 $$
where $M$ is locally free as a sheaf on $F$. Then $2 \chi(M) \geq e + \nu(M)$ where 
$$ \nu(M):= \chi(\shom_{\bara}(M,M)) - \chi(\shom_{\bara}(M',M')) 
- \chi(\sext^1_{\bara}(M,M)) + \chi(\sext^1_{\bara}(M,M')).$$
\end{prop}
\textbf{Proof.} We use the long exact sequences
\begin{align*} 
0 \lm \shom_{\bara}(M,M) \lm M \lm \shom_{\bara}(M',M) \lm \sext^1_{\bara}(M,M) \lm 0 \\
0 \lm \shom_{\bara}(M,M') \lm M' \lm \shom_{\bara}(M',M') \lm \sext^1_{\bara}(M,M') \lm 0.
\end{align*}
From proposition~\ref{pnoExt}ii), we obtain the inequality
\begin{align*}
0 & \leq \chi(\shom_{\bara}(M',M))) - \chi(\shom_{\bara}(M,M')) \\
  & = \chi(M) - \chi(\shom_{\bara}(M,M)) + \chi(\sext^1_{\bara}(M,M)) - \chi(M') + \chi(\shom_{\bara}(M',M')) - \chi(\sext^1_{\bara}(M,M'))
\end{align*}
The proposition now follows on noting that $\chi(M') + \chi(M) = \chi(\bara) = e$ by continuity of Euler characteristic and theorem~\ref{tgeneric}. 

\noindent
\textbf{Remark:} As will be seen in lemma~\ref{ltannu}, one can compute $\nu(M)$ locally.

\section{Hereditary fibres}  \label{shered}  

We return to the study of our degree $e$ ruled order $A$ on the ruled surface $f:Z \lm C$. In this section, we examine goodness of fibres $F$ of $f:Z \lm C$ which meet the ramification locus transversally. This corresponds to the case where $\bara:=A|_F$ is an hereditary order. 

We start with a lemma we will use more generally.
\begin{lemma} \label{ldijdji}  
Let $\bara$ be a degree $e$ hereditary order on $F=\PP^1$ which is totally ramified at some point $p$ (that is, the ramification index of $\bara$ at $p$ is $e$). Suppose that $\bara$ has a complete set of $e$ idempotents  and that the corresponding Peirce decomposition has the form
$$\bara = \bigoplus_{i,j = 1}^{e} \calo_{\PP^1}(p_{ij}) .$$
Then for $i \neq j$ we have $\deg p_{ij} + \deg p_{ji} < 0$. 
\end{lemma}
\textbf{Proof.} Let $d_{ij} = \deg p_{ij}$. Note that closure under multiplication ensures that $d_{ij} + d_{ji} \leq 0$. We wish to derive a contradiction from the assumption $d_{ij} + d_{ji} = 0$. Let $\eps\in H^0(\bara)$ be the idempotent which has a 1 in the $i$ and $j$-th diagonal entry and 0s elsewhere. Now $\eps \bara \eps$ is the Azumaya algebra $\send (\calo \oplus \calo(d_{ij}))$ so $\eps (\bara \otimes_F k(p)) \eps \simeq k^{2 \times 2}$. 
If $J$ is the Jacobson radical of the ring $\bara \otimes_F k(p)$, then since $\bara$ is totally ramified at $p$, we know that $\bara \otimes_F k(p)/J \simeq k^e$. However, the nilpotent ideal $\eps J \eps\triangleleft \eps (\bara\otimes_F k(p)) \eps$ must be 0. This contradicts the fact that $\bara \otimes_F k(p)/J$ is commutative. The lemma is proved. 

\begin{thm}  \label{theredfibres}  
Let $A$ be a degree $e$ ruled order such that $c_2(A)$ is minimal in its Morita equivalence class. Then for every fibre $F$ which meets the ramification locus transversally, we have $H^1(\bara) = 0$ so the results of theorem~\ref{tgeneric} apply. In particular, $F$ is good.
\end{thm}
\textbf{Proof.} From lemma~\ref{lheredPeirce}, there is a Peirce decomposition 
$$\bara = \bigoplus_{i,j = 1}^{e} \calo_{\PP^1}(p_{ij}) $$
for some divisors $p_{ij} \in \Div \PP^1$ of degree $d_{ij}$. Note that $p_{ii} = 0$ for all $i$. Let $P_j$ be the locally projective $\bara$-module corresponding to the $j$-th column $\oplus_i \calo(p_{ij})$. Note that $\shom_{\bara}(P_i,P_j) = \calo(p_{ij})$.

We wish first to show that $d_{ij} <0$ for $i \neq j$. To this end, let $A'$ be the Morita transform of $A$ associated to the quotient $P_j$. We use proposition~\ref{pMorita} with $M= P_j, M'= \oplus_{l\neq j} P_l$ both locally projective to see
$$ 0 \leq c_2(A') - c_2(A) = \sum_{l\neq j} d_{lj} - \sum_{l\neq j} d_{jl}.$$
Adding $d_{jj}$ to both sums shows that $\sum_{l} d_{lj} \geq \sum_l d_{jl}$. We can also apply a Morita transform to $A$ associated to the quotient $\oplus_{l\neq j} P_l$ to obtain the reverse inequality so 
$$\sum_l d_{lj} = \sum_l d_{jl}.$$

We argue by contradiction and assume without loss of generality that $d_{12} \geq 0$ so by closure under multiplication we have $d_{l2}\geq d_{l1}, d_{1l} \geq d_{2l}$ for any $l$. Then 
$$\sum_l d_{l2} \geq \sum_l d_{l1} = \sum_l d_{1l} \geq \sum_l d_{2l} = \sum_l d_{l2} .$$
Hence equality must hold throughout, which in turn implies that $d_{l2}= d_{l1}, d_{1l}= d_{2l}$ for every $l$. In particular, $d_{12} = 0 = d_{21}$. However, $\bara$ is totally ramified at two points (since $A$ has ramification index $e$ along the bisection $D$) so the lemma gives a    contradiction. We have thus proved that $d_{ij} <0$ for $i \neq j$ so $H^0(\bara) = e$. Furthermore, $A$ is flat over $C$ so $\chi(\bara) = e$ showing $H^1(\bara) = 0$. This completes the proof of the theorem.

\section{Minimal hereditary orders containing $\bara$}  \label{sheredover}  

Our approach to analysing the non-hereditary fibres $\bara$ of $A$ is to embed them in an hereditary order $\bara_1$ which is close to $\bara$ and then bound the possibilities for $\bara_1$. Our setup is as follows. As usual, let $A$ be a degree $e$ ruled order on the ruled surface $f:Z \lm C$ with minimal $c_2$ in its Morita equivalence class. Pick a closed fibre $F$ which intersects the ramification divisor in exactly one point $p$ (necessarily then with multiplicity 2). Suppose there is an hereditary order $\bara_1$ containing $\bara:= A|_F$ which is totally ramified at $p$. By lemma~\ref{lheredPeirce}, there is a Peirce decomposition $\bara_1 = \oplus \calo(p_{ij})$ and the indecomposable summands $P_1,\ldots, P_e$ of $\bara_1$ correspond to columns of this matrix form.

The following result is standard though we do not know a suitable reference for it. 
\begin{prop} \label{pheredproj}  
Let $P,P'$ be locally projective $\bara_1$-modules of rank $e$. 
\begin{enumerate}
\item If $d,d'$ are the degrees of the invertible $\calo_F$-modules $\shom_{\bara_1}(P,P'),\shom_{\bara_1}(P',P)$, then $d+d'=0$ if $P \otimes_F \hat{\calo}_{F,p} \simeq P' \otimes_F \hat{\calo}_{F,p}$ and is $-1$ otherwise.
\item If $\chi(P) = \chi(P')$, then $P \simeq P'$.
\end{enumerate}
\end{prop}
\textbf{Proof.} We prove part~i) first by analysing the algebra $\send_{\bara_1}(P \oplus P')$.  Since $\bara_1$ is totally ramified at $p$, we have 
$$ \bara_{1p} \simeq 
\begin{pmatrix}
k[[z]] & \cdots & \cdots & k[[z]] \\
(z)      & k[[z]] &       & \vdots \\
\vdots   & \ddots & \ddots  & \vdots \\
(z)      & \cdots & (z)     & k[[z]]
\end{pmatrix}
$$

Suppose $P \otimes_F \hat{\calo}_{F,p},P' \otimes_F \hat{\calo}_{F,p}$ are isomorphic to the $i$-th and $j$-column of the matrix form above. If $i = j$ then $\send_{\bara_1}(P \oplus P') \otimes_F \hat{\calo}_{F,p}$ is the full matrix algebra in $\hat{\calo}_{F,p}$. If $i\neq j$ then let $\eps\in \bara_{1p}$ be the idempotent with 1's in the $i$-th and $j$-th diagonal entry and 0's elsewhere. Then $\send_{\bara_1}(P \oplus P') \otimes_F \hat{\calo}_{F,p}$ is the hereditary algebra $\eps \bara_{1p} \eps$ with discriminant ideal $(z)$. 

Looking globally we see that 
$$ \send_{\bara_1}(P \oplus P') \simeq
\begin{pmatrix}
\calo_F & \calo_F(q') \\ \calo_F(q) & \calo_F
\end{pmatrix}
$$
for some divisors $q,q'$ of degrees $d,d'$. Conjugating the above matrix if necessary, we may assume that $q=dp$. Furthermore, $\send_{\bara_1}(P \oplus P')$ is Azumaya away from $p$ since the same is true of $\bara_1$. Hence we may also assume $q'=d'p$. The local computations in the previous paragraph now show that $d+d'=0$ if $i=j$ and $d+d'=-1$ otherwise.

We now prove part~ii). From part~i), we may assume without loss of generality that $d\geq 0$, so $\Hom_{\bara_1}(P,P') \neq 0$. Then any non-zero morphism $P \lm P'$ must be an isomorphism, since both are 1-critical and have the same Euler characteristic. 
\vspace{2mm}

We now consider Morita transforms arising from the following setup. Let $M_1$ be a locally projective quotient of $\bara_1$ and $M$ be the image of $\bara$ under $\bara_1 \lm M_1$. We complete the following commutative diagram with exact rows.
\begin{equation*}
\diagram
0 \rto & M' \dto\rto & \bara \dto\rto & M \dto\rto & 0 \\
0 \rto & M'_1 \rto & \bara_1 \rto & M_1 \rto & 0 
\enddiagram
\label{eheredover}
\end{equation*}
In the next lemma we use the number
$$\nu(M) = \chi(\shom_{\bara}(M,M)) - \chi(\shom_{\bara}(M',M')) 
- \chi(\sext^1_{\bara}(M,M)) + \chi(\sext^1_{\bara}(M,M'))$$
introduced in proposition~\ref{pminc2ruled}.
\begin{lemma}  \label{lbdhered}  
We assume the above setup, 
\begin{enumerate}
\item Suppose that when $M_1=P_i$, we have $\nu(M) \geq 2-e$. Then $\chi(P_i) \geq 1$.
\item  Suppose that $\chi(P_i) \geq 1$ for all $i$. Then 
$$ \bara_1 \simeq 
\begin{pmatrix}
\calo & \cdots & \cdots & \calo \\
\calo(-1)      & \calo &       & \vdots \\
\vdots   & \ddots & \ddots  & \vdots \\
\calo(-1)      & \cdots & \calo(-1)     & \calo
\end{pmatrix}.$$
\end{enumerate} 
\end{lemma}
\textbf{Proof.} To see i), note that when $M_1=P_i$, the assumption on $\nu$ and proposition~\ref{pminc2ruled}, show that $\chi(P_i) \geq \chi(M) \geq 1$.

To prove ii), we first show that the components $\calo(p_{ij})$ of the Peirce decomposition for $\bara_1$ satisfy
\begin{claim}  \label{cdij} 
Let $d_{ij}:=\deg p_{ij}$. Then for distinct $i,j$ we have $\{d_{ij},d_{ji}\} = \{-1,0\}$.
\end{claim}
\textbf{Proof.} Complete locally at $p$, the $P_i$ are the $e$ distinct indecomposable projective $\bara_{1p}$-modules, so by proposition~\ref{pheredproj}i), it suffices to show that $d_{ij} \leq 0$. Now 
if $d_{ij} \geq 1$, then $\chi(P_j) \geq \chi(P_i) + e$ so it in turn suffices to show that for all $i$ we have $\chi(P_i) \leq e$.  Proposition~\ref{pheredproj} ensures that $\chi(P_1),\ldots,\chi(P_e)$ are distinct so the assumption on the $\chi(P_i)$ forces
$$ \chi(\bara_1) = \sum_{i=1}^e \chi(P_i) \geq 1 + 2 + \ldots + e =  \frac{1}{2}e(e+1)$$
with equality if and only if $\{\chi(P_1),\ldots,\chi(P_e)\} = \{1,\ldots,e\}$. But we may also compute $\chi(\bara_1)$ as follows. Since $\bara_1$ is totally ramified at $p$, we know from the local structure theory of hereditary orders, that we may embed $\bara_1$ in a maximal order $\bara_2$ so that $\chi(\bara_2/\bara_1) = \frac{1}{2}e(e-1)$. Also, Tsen's theorem tells us that $\bara_2$ is trivial Azumaya so $\chi(\bara_2) = e^2$ giving $\chi(\bara_1) = \frac{1}{2}e(e+1)$. This proves the claim.

Returning to the proof of the lemma, we wish to show that after conjugating $\bara_1$ by some permutation matrix, we have $d_{ij} = 0$ for $i \leq j$ and $d_{ij} = -1$ for $i > j$. To find this permutation matrix, note that we can partially order the set of columns $\{P_i\}^e_{i=1}$ by $P_i \leq P_j$ if $d_{ij} = 0$. Indeed, the relation is reflexive since $d_{ii} = 0$, anti-symmetric by the claim and transitive by closure under multiplication. Furthermore, the claim shows that the order on $\{P_i\}$ is in fact a total order. Conjugating by a permutation matrix, we may thus assume that $P_i \leq P_j$ if and only if $i \leq j$ which proves the lemma.

\section{Rational filtrations and goodness}  \label{sratgood}  

In this section, we introduce the notion of rational filtrations to give a criterion for goodness of a fibre.

We let $\bara$ be a coherent sheaf of algebras on $F = \PP^1$ which is locally free of rank $e^2$ as a sheaf. Suppose further that $k(\bara):= \bara \otimes_{F} k(F)$ satisfies
\begin{ass}  \label{arankesimples} 
All simple modules have dimension $e$ over $k(F)$.
\end{ass}
In particular, any $\bara$-module which is rank $e$ torsion-free is 1-critical. Note that any order in the full matrix algebra $M_e(k(F))$ satisfies the hypothesis.

We first make the following 
\begin{defn}  \label{dratmod}  
An $\bara$-module $N$ is said to be {\em rational} if it is torsion-free of rank $e$ and has cohomology $h^0(N) = 1, h^1(N) = 0$. Equivalently, we have the isomorphism of sheaves $N \simeq \calo \oplus \calo(-1)^{\oplus e-1}$. A filtration 
$$ 0 \leq N^1 \leq N^2 \leq \ldots \leq N^r = N $$
of $N$ is said to be {\em rational} if all the factors $N^{i+1}/N^i$ are rational or zero, in which case we also say that $N$ is {\em rationally filtered}.
\end{defn}

\begin{prop}  \label{pratfilt}  
Any rationally filtered module $N$ satisfies $H^1(N) = 0, \rk N = eh^0(N)$. 
\end{prop}
\textbf{Proof.} As a sheaf on $F$, $N$ is just the direct sum of the factors in a rational filtration. The proposition follows. 

\vspace{2mm}

Rationally filtered modules enjoy the following nice properties.
\begin{lemma} \label{lmorphismrat}  
Let $\phi:N_1 \lm N_2$ be a non-zero morphism of $\bara$-modules where $N_1$ is rationally filtered. Then 
\begin{enumerate}
\item $H^1(\im \phi) = 0$. 
\item if $N_2$ is rational, then $\phi$ is surjective whenever it is non-zero.
\item if $N_2$ is rationally filtered then so are $\ker \phi$ and $\coker \phi$. 
\end{enumerate} 
\end{lemma}
\textbf{Proof.} Part i) follows from proposition~\ref{pratfilt} and the long exact sequence in cohomology. To prove part ii), suppose now that $N_2$ is rational and $\phi$ is non-zero so that $\im \phi$ has rank $e$ too. Grothendieck's splitting theorem and part i) ensure $\im \phi \simeq \oplus_{i=1}^e \calo(d_i)$ where the $d_i \geq -1$. If $\phi$ is not surjective then we must have $\im \phi \simeq \calo(-1)^{\oplus e}$. We derive a contradiction as follows. Let $N'_1$ be a rational submodule of $N_1$. The restricted map $\phi':N'_1 \lm N_2$ is not injective since the summand $\calo$ maps to zero. Hence $\phi'=0$. Continuing inductively up the filtration we see that $\phi = 0$. 

Finally we prove part iii) and so assume there are rational filtrations
\begin{align*}
0 < N^1_1 < N^2_1 < \ldots < N^r_1 = N_1 \\
0 < N^1_2 < N^2_2 < \ldots < N^s_2 = N_2 
\end{align*}
We argue by induction on $s$ and assume first that $s=1$. We may as well assume that $\phi$ is non-zero and so can pick $i$ maximal such that $\phi(N^i_1) = 0$. The induced map $N^{i+1}_1/N^i_1 \lm N_2$ is thus an isomorphism by part ii), and its inverse yields a splitting of the map $N_1/N^i_1 \lm N_2$. Hence $\ker \phi$ is an extension of $N_1/N^{i+1}_1$ by $N^i_1$ so is rationally filtered. This completes the case $s=1$. 

For general $s$, we may assume by induction that $\phi$ does not factor through $N^{s-1}_2$ in which case the composite $N_1 \lm N_2 \lm N_2/N_2^{s-1}$ is surjective by part~ii). We may thus consider the following commutative diagram with exact rows
$$\diagram
0 \rto &  K \rto \dto^{\phi'} & N_1  \rto\dto^{\phi} & N_2/N_2^{s-1}  \rto\ddouble & 0 \\ 
0 \rto &  N_2^{s-1} \rto & N_2  \rto^(.35){\pi} & N_2/N_2^{s-1}  \rto & 0 \\ 
\enddiagram$$
where $\pi$ is the natural quotient map and $K$ is the appropriate kernel. The $s=1$ case shows that $K$ is also rationally filtered. Now $\ker\phi = \ker \phi', \coker \phi = \coker \phi'$ and both are rationally filtered by the inductive hypothesis. The lemma is proved.
\vspace{2mm} 

\begin{thm}  \label{tratgood}  
Let $A$ be a degree $e$ ruled order on a ruled surface $Z$ and $F$ be a ruling of $Z$ such that $\bara := A|_F$ satisfies hypothesis~\ref{arankesimples} and is rationally filtered as an $\bara$-module. Then $F$ is a good fibre.
\end{thm}
\textbf{Proof.} By proposition~\ref{pflatcrit} and our hypothesis~\ref{arankesimples}, we need only show that for any algebra quotient map $H^0(\bara) \lm k$ we have $\bara \otimes_{H^0(\bara)} k \simeq \calo \oplus \calo(-1)^{\oplus e-1}$. Let $J \triangleleft H^0(\bara)$ be the Jacobson radical. 
Now $H^1(\bara)=0$ so cohomology commutes with base change and we see that $H^0(\bara)$ is commutative. Hence $H^0(\bara)/J \simeq \prod_{l=1}^r k\eps_l$ for some idempotents $\eps_l$. We need to show that the direct summand $(\bara/\bara J) \eps_l$ of $\bara/\bara J$ is rational. First note that $\bara/\bara J$ is rationally filtered by lemma~\ref{lmorphismrat}iii) applied to $\bara \otimes_k J \lm \bara$, so the same is true of $(\bara/\bara J) \eps_l$. Also, $(\bara/\bara J) \eps_l\neq 0$ since it contains $\eps_l$. It suffices now to check its rank. Now $H^1(\bara J) = 0$ since $\bara J$ is a quotient of a direct sum of $\dim J$ copies of $\bara$. Hence
$$ \rk (\bara/\bara J) = eh^0(\bara/\bara J) = eh^0(\bara) - eh^0(\bara J) \leq eh^0(\bara) - e\dim J = e\dim H^0(\bara)/J = er .$$
It follows that the summand $(\bara/\bara J) \eps_l$ must have rank $e$ and so be rational. This proves the theorem.

\section{Fibres tangential to $D$}  \label{sfibretan} 

We assume as usual that $A$ is a degree $e$ ruled order with minimal $c_2$ in its Morita equivalence class. In this section, we consider a fibre $F$ which is tangential to the ramification curve $D$ and let $p$ be the intersection point $F \cap D$. We use a subscript $p$ to denote the completion at $p$, so for example $\bara_p := \bara \otimes_F \hat{\calo}_{F,p}$. Theorem~\ref{tlocalterm} can be used to determine $\bara_p$. 

We will pick an hereditary over-order $\bara_1$ of $\bara$ so that the inclusion $\bara_p \subset \bara_{1p}$ is given by  
$$\begin{pmatrix}
k[[z]] & \cdots & \cdots & k[[z]] \\
(z^2)      & k[[z]] &       & \vdots \\
\vdots   & \ddots & \ddots  & \vdots \\
(z^2)      & \cdots & (z^2)     & k[[z]]
\end{pmatrix} \subset 
\begin{pmatrix}
k[[z]] & \cdots & \cdots & k[[z]] \\
(z)      & k[[z]] &       & \vdots \\
\vdots   & \ddots & \ddots  & \vdots \\
(z)      & \cdots & (z)     & k[[z]]
\end{pmatrix}$$
where we have used the isomorphism $\calo_{F,p} \simeq k[[z]]$. Note that $\bara_1$ is an hereditary order which is totally ramified at the point $p$, so we are in a good position to use the results of sections~\ref{sheredover} and \ref{sratgood}. 

Since $\bara_p$ has a complete set of idempotents, the rank $e$ torsion-free $\bara_p$-modules are easily computed to be (up to isomorphism), those of the form 
$$  P_{ij} := 
\begin{pmatrix}
I_1 \\ \vdots \\ I_e
\end{pmatrix}, \ \ i,j \in \{ 1,\ldots, e\}$$
where $I_1 = \cdots = I_i = k[[z]], I_{i+1} = \cdots = I_{j} = (z), I_{j+1} = \cdots = I_e = (z^2)$. The projective $\bara_p$-modules are those of the form $P_{ii}$ which correspond to the columns of the matrix form. The columns of $\bara_{1p}$ are those of the form $P_{ie}$. 

First note that we have the following non-split exact sequences
\begin{eqnarray}
0 \lm P_{rs} \lm P_{js} \oplus P_{ir} \lm P_{ij} \lm 0, \text{for}\ 
i \leq r < j \leq s \label{e1} \\
0 \lm P_{rs} \lm P_{is} \oplus P_{rj} \lm P_{ij} \lm 0, \text{for}\ 
r < i \leq s < j \label{e2}
\end{eqnarray}
We can now classify the indecomposable torsion-free $\bara_p$-modules.
\begin{prop}  \label{pindec}  
\begin{enumerate}
\item $\Ext_{\bara_p}^1(P_{ij},P_{rs}) = k$ if $i \leq r < j \leq s$ or $r < i \leq s < j$ and is zero otherwise.
\item Every indecomposable torsion-free $\bara_p$-module is isomorphic to $P_{ij}$ for some $i,j$.
\end{enumerate}
\end{prop}
\textbf{Proof.} i) Equation~(\ref{e1}) gives in particular a partial projective resolution
$$  0 \lm P_{ij}  \lm P_{ii} \oplus P_{jj} \lm P_{ij}  \lm 0 $$
from which we can calculate the ext groups.

ii) Let $N$ be a torsion-free $\bara_p$-module. We show it is the direct sum of rank $e$ torsion-free modules by induction on rank. Pick a surjection $N \lm P_{ij}$ with $j-i$ minimal. Such surjections exist since $\bara_p$ is an order in a full matrix algebra over $\calo_{F,p}$. We may assume the surjection is not split, otherwise we are done by induction. Consider the corresponding exact sequence
$$ E:0 \lm K \lm N \lm P_{ij} \lm 0 $$
which we view as a non-zero extension $E \in \Ext^1_{\bara_p}(P_{ij},K)$. By induction, $K = \oplus_l P_{r_ls_l}$. Pick some $l$ such that for $r=r_l,s=s_l$, the image of $E$ in $\Ext^1_{\bara_p}(P_{ij},P_{rs})$ is non-zero. We obtain thus a morphism of extensions
$$\diagram
0 \rto & K \rto\dto & N \rto\dto & P_{ij} \rto\ddouble & 0 \\
0 \rto & P_{rs} \rto & N' \rto & P_{ij} \rto & 0  
\enddiagram$$
Now from part i) we know $\Ext_{\bara_p}^1(P_{ij},P_{rs}) = k$ and either $i \leq r < j \leq s$ or $r < i \leq s < j$. In the first case, we have as in equation~(\ref{e1}) that $N' \simeq P_{js} \oplus P_{ir}$. This is a contradiction since we then obtain a surjection $N \lm N' \lm P_{ir}$ and $r-i < j-i$. In the second case, we have as in equation~(\ref{e2}) that $N' \simeq P_{is} \oplus P_{rj}$ and again obtain a similar contradiction.

\vspace{2mm}

We omit the elementary proof of the next presumably well-known result.

\begin{lemma}  \label{lsameker}
Let $R$ be a noetherian ring with Jacobson radical $J$ and $P$ be a finitely generated projective $R$-module. Suppose that $R/J$ is artinian. Given any submodule $P'<P$ and surjective homomorphism $\phi: P \lm P/P'$ we have $\ker \phi \simeq P'$. 
\end{lemma}

The possibilities for the exact sequence $0 \lm M' \lm \bara \lm M \lm 0$ complete locally at $p$ are given in the next result. We introduce the following notation. If $\eps_1,\ldots, \eps_e \in \bara_p$ are the standard diagonal idempotents with exactly one non-zero entry, then the {\em dimension vector} of any $\bara_p$-module $N$ is defined to be $(\dim_k \eps_1 N, \ldots, \dim_k \eps_e N)$. For example, the dimension vector of $P_{ij}/(\text{rad}\ \bara_p)P_{ij}$ is $(0,\ldots, 0,1,0,\ldots,0,1,0 \ldots,0)$ where the 1s occur in the $i$-th and $j$-th co-ordinate. 
\begin{prop}  \label{pquotAbar} 
Let $P$ be a torsion-free $\bara_p$-module and $\phi: \bara_p \lm P$ be a surjective morphism. Then there is a partition of $\{1,\ldots,e\}$ into 4 disjoint subsets
$$ \{i_1,\ldots, i_m\}, \{j_1,\ldots, j_m\}, \Lambda, \Lambda'$$ 
such that 
\begin{enumerate}
\item $P \simeq \bigoplus_{\lambda \in \Lambda} P_{\lambda\lambda} \oplus \bigoplus_{l} P_{i_lj_l}$.
\item $\ker \phi \simeq \bigoplus_{\lambda \in \Lambda'} P_{\lambda\lambda} \oplus \bigoplus_{l} P_{i_lj_l}$.
\end{enumerate}
\end{prop} 
\textbf{Proof.} Part i) follows from proposition~\ref{pindec}ii) and the fact that $\bara_p/\text{rad} \bara_p$ has dimension vector $(1,1,\ldots,1)$. We now prove part ii). By lemma~\ref{lsameker}, we may assume that $\phi$ is a direct sum of surjections of the form $P_{i_li_l} \oplus P_{j_lj_l} \lm P_{i_lj_l}$, $P_{\lambda\lambda} \lm P_{\lambda\lambda}$ for $\lambda \in \Lambda$ and $P_{\lambda\lambda} \lm 0$ for $\lambda \in \Lambda'$. The proposition follows now from equation~(\ref{e1}). 

\vspace{2mm}

We wish to use lemma~\ref{lbdhered} to show $H^1(\bara_1) = 0$ first. We need to calculate
$$\nu(M) = \chi(\shom_{\bara}(M,M)) - \chi(\shom_{\bara}(M',M')) 
- \chi(\sext^1_{\bara}(M,M)) + \chi(\sext^1_{\bara}(M,M'))$$
which is possible since it depends only on local data.

\begin{lemma}  \label{ltannu}  
Consider an exact sequence of $\bara$-modules
$$ 0 \lm M' \lm \bara \lm M \lm 0$$
with $M$ torsion-free. Then
\begin{enumerate}
\item $\chi(\sext^1_{\bara}(M,M)) = \chi(\sext^1_{\bara}(M,M'))$, and
\item $\chi(\shom_{\bara}(M,M)) - \chi(\shom_{\bara}(M',M'))= \frac{1}{e}(\rk M - \rk M')$. 
\end{enumerate}
In particular, $\nu(M) = \frac{1}{e}(\rk M - \rk M')$. 
\end{lemma}
\textbf{Proof.} To prove i), note that $M$ is locally projective as an $\bara$-module away from $p$, so the ext groups can be computed after completing at $p$. Now proposition~\ref{pindec}i) shows $\Ext^1_{\bara_p}(P_{ij},P_{rr}) = 0$ for all $i,j,r$ while proposition~\ref{pquotAbar} shows that $M_p$ and $M'_p$ only differ by projectives so part i) holds. 

We now verify ii). Suppose that $M$ has rank $re$. We can embed $\shom_{\bara} (M,M)$ in a maximal order $B$ where locally $B_p \simeq \calo^{r \times r}_{F,p}$. Now $\shom_{\bara}(M,M)$ is maximal everywhere except at $p$ so 
$$\chi(\shom_{\bara}(M,M)) = \chi(B) - \text{colength}\ \shom_{\bara_p}(M_p,M_p) = r^2 - \text{colength}\ \shom_{\bara_p}(M_p,M_p)$$
where the colength is computed as a subsheaf of $B_p$. To compute the colength term, we use proposition~\ref{pquotAbar} to write 
$$ M_p \simeq \bigoplus_{\lambda \in \Lambda} P_{\lambda\lambda} \oplus \bigoplus_{l=1}^m P_{i_lj_l}, M'_p \simeq \bigoplus_{\lambda \in \Lambda'} P_{\lambda\lambda} \oplus \bigoplus_{l=1}^m P_{i_lj_l} $$
for a partition $\{i_1,\ldots, i_m\}, \{j_1,\ldots, j_m\}, \Lambda, \Lambda'$ of $\{1, \ldots e\}$. The direct sum decomposition for $M_p$ gives a Peirce decomposition for $\shom_{\bara_p}(M_p,M_p)$. If we use an isomorphism $\calo_{F,p} \simeq k[[z]]$, we can write some of these Peirce components as follows. 
\begin{enumerate}
\item $\shom_{\bara_p}(P_{\lambda\lambda},P_{\rho\rho}) \oplus \shom_{\bara_p}(P_{\rho\rho},P_{\lambda\lambda}) \simeq k[[z]] \oplus (z^2)$ for $\lambda \neq \rho$. 
\item $\shom_{\bara_p}(P_{\lambda\lambda},P_{ij}) \oplus \shom_{\bara_p}(P_{ij},P_{\lambda\lambda}) \simeq (z) \oplus (z)$ if $i< \lambda \leq j$ and $\simeq k[[z]] \oplus (z^2)$ otherwise.
\end{enumerate}
In either case, their colength (in some appropriate copy of $\calo^2 \subset \calo^{r \times r}$) is 2. We thus find 
$$\chi(\shom_{\bara}(M,M))  = r^2 - \text{colength}\ \shom_{\bara_p}(M_p,M_p) = r - \text{colength}\ \send_{\bara_p} P' +m^2 - m$$
where $P' = P_{i_1j_1} \oplus \ldots \oplus P_{i_mj_m}$. Similarly, if $M'$ has rank $r'e$ we find
$$\chi(\shom_{\bara}(M',M')) = r' - \text{colength}\ \send_{\bara_p} P' +m^2 - m$$
which completes the proof of the lemma.
\vspace{2mm}

\begin{thm}  \label{theredtan}  
$$ \bara_1 \simeq 
\begin{pmatrix}
\calo & \cdots & \cdots & \calo \\
\calo(-1)      & \calo &       & \vdots \\
\vdots   & \ddots & \ddots  & \vdots \\
\calo(-1)      & \cdots & \calo(-1)     & \calo
\end{pmatrix}.$$
\end{thm}
\textbf{Proof.} This follows from lemma~\ref{lbdhered} and the $\nu$ computations in lemma~\ref{ltannu}. 
\vspace{2mm}

We wish now to show that $\bara$ is rationally filtered. To this end, let $P_i\simeq \calo^{\oplus i} \oplus \calo(-1)^{\oplus e-i}$ be the indecomposable summand of $\bara_1$ corresponding to the $i$-th column of the matrix form in theorem~\ref{theredtan}. We work a little more generally and consider an arbitrary order $B$ in $\bara_1$.  We define 
$$ B^i := B \cap P_{>e-i} , \ \ 
S^i := \text{coker}(B^i \hookrightarrow P_{>e-i}) \label{filtration_on_B}$$
where $P_{>e-i} = P_{e-i+1} \oplus \ldots \oplus P_e$. 
We thus obtain a filtration
$$ 0 < B^1 < \ldots < B^e := B.$$
When $B = \bara$, we wish to apply proposition~\ref{pminc2ruled} in the following situation. Let $\phi: \bara \hookrightarrow \bara_1 \lm P_{\leq e-i}:= P_1 \oplus \ldots \oplus P_{e-i}$ be the composite of the natural inclusion with the natural projection. Then we can apply a Morita transform with respect to $M = \im \phi$ and note that $M'=\ker \phi = \bara^i$. 
\begin{lemma}  \label{ltanfilt}  
$\chi(\bara^i) = i$ so $\chi(\bara^{i+1}/\bara^i) = 1$. 
\end{lemma}
\textbf{Proof.}  Using the notation above, we see from lemma~\ref{ltannu} and proposition~\ref{pminc2ruled}, that $\chi(M)\geq \frac{1}{e}\rk M$. Hence $\chi(\bara^i) \leq e - \frac{1}{e}\rk M = \frac{1}{e}\rk \bara^i = i$. We proceed to prove the reverse inequality. Now $S^i$ is both a quotient of $P_{>e-i} \otimes_F k(p)$, which has dimension vector $(i,i,\ldots, i)$, and a submodule of $\bara_1/\bara$, which has dimension vector $(0,1,2,\ldots,e-1)$. Hence 
$$\chi(S^i) \leq 0 + 1 + 2 + \ldots + i + i + \ldots + i$$ 
where the number of $i$'s on the right hand side is $e-i$. 
Also $\chi(P_{>e-i}) = e + (e-1) + \ldots + (e-i+1)$ so 
$$\chi(\bara^i) = \chi(P_{>e-i}) - \chi(S^i) \geq i.$$
This completes the proof of the lemma. 
\vspace{2mm}

\begin{thm} \label{ttannoH1}  
Let $B \subset \bara_1$ be an order with $\chi(B^i) = i$. Then 
the filtration $\{B^i\}_{i=1}^e$ is rational. In particular, $\bara$ is rationally filtered and $F$ is a good fibre. 
\end{thm}
\textbf{Proof.}  Note that $\chi(B^{i+1}/B^i) = 1$ so in particular, the injection $B/B^{e-1} \hookrightarrow P_1$ must be an isomorphism. We proceed by downward induction on $i$. 
Let $i$ be maximal such that $B^{i+1}/B^i\simeq (B^{i+1} + P_{>e-i})/P_{>e-i}$ is not rational. Then we must have $h^0(B^{i+1}/B^i) \geq 2$. Also, $P_{>e-i}$ is a direct sum of line bundles isomorphic to $\calo$ and $\calo(-1)$ so we may lift sections to find $s \in H^0(B^{i+1} + P_{>e-i})$ such that in the matrix form of theorem~\ref{theredtan}, $s$ is zero in the first $e-i-1$ columns and in the $(e-i)$-th column is i) zero on or below the main diagonal and ii) non-zero strictly above the main diagonal. An elementary matrix computation shows that right multiplication by $s$ induces a non-zero morphism $\s:B/B^{i+1} \lm (B^{i+1} + P_{>e-i})/P_{>e-i}$. Now $B/B^{i+1}$ is rationally filtered so lemma~\ref{lmorphismrat} shows that $H^1(\im \s) = 0$. Also, the cokernel of $\im \s \lm (B^{i+1} + P_{>e-i})/P_{>e-i}$ has finite length, so $H^1(B^{i+1}/B^i) = 0$ too. This gives a contradiction as then $\chi(B^{i+1}/B^i) = h^0(B^{i+1}/B^i)> 1$. Thus $\{B^i\}$ is a rational filtration and the rest of the theorem follows from lemma~\ref{ltanfilt} and theorem~\ref{tratgood}. 
\vspace{2mm}

\section{Fibres intersecting a node of $D$\label{sfibrenode}}

Let $A$ be a ruled order of degree $e$ on $f : Z \to C$ with
minimal $c_2$ in its Morita equivalence class. Let $F$ be a fibre of $f$ which
intersects a node $p$ of the ramification divisor $D$. Following the notation of the previous section, a subscript $p$ will denote completion at $p$ and we fix an isomorphism
$\mathcal{O}_{F, p} \simeq k [[z]]$. Let $\bara  = A \otimes_Z \mathcal{O}_F$ and $k ( \bara _p)$ be the central
simple $k((z))$-algebra $\bara _p \otimes_F k (F)$.

By theorem \ref{tlocalterm}, we have the
following isomorphism
\begin{eqnarray}
  \bara _p & \simeq & \frac{k [[z]] \left\langle x, y \right\rangle}{(x^e -
  z, y^e - z, x y - \zeta y x)}  \label{Appresentation}
\end{eqnarray}
where $\zeta$ is a primitive $e$-th root of unity. The Jacobson radical $J$ of
$\bara _p$ is generated by $x$ and $y$, so $\bara _p / J \simeq k$ and
$\bara _p$ is a local ring. Following section \ref{sheredover}, we wish to study
$\bara _p$ by embedding it in an hereditary order. To do this we first
determine the torsion-free $\bara _p$-modules of rank $e$ up to isomorphism.

\begin{lemma}
  \label{I-adicfiltrationonP}Let $I \subseteq J$ be an ideal of $\bara _p$
  such that $z \bara _p \subseteq I^e$, and $Q$ be a torsion-free
  $\bara _p$-module of rank $e$. Then the $I$-adic filtration on $Q$ has
  $1$-dimensional quotients. In particular, the $I$-adic and $J$-adic
  filtrations on $Q$ are equal. 
\end{lemma}
\textbf{Proof.}  
  Note that $I \, Q$ is torsion-free of rank $e$, so it suffices by induction to show
  that $Q / I \, Q \simeq k$. Since $z \bara _p \subseteq I^e$ we have a
  surjection $Q / z \, Q \to Q / I^e Q$, so $\dim_k (Q / I^e Q)
  \leqslant \dim_k (Q / z \, Q) = e$. Nakayama's lemma gives $I^k Q / I^{k +
  1} Q \neq 0$, and since $\dim_k (Q / I^e Q) = \sum_{k = 0}^{e - 1} \dim_k
  (I^k Q / I^{k + 1} Q)$, we have $Q / I \, Q \simeq k$. From this it follows
  that $I \, Q = J \, Q$, so the $I$-adic and $J$-adic filtrations on $Q$ are
  equal. 
\vspace{2mm}

\begin{prop}
  \label{defQi}Let $Q_i$ denote the quotient $\bara _p / \bara _p \xi_i$
  where $\xi_i = x - \zeta^{i - (e - 1) / 2} y$. Then $Q_i$ is a torsion-free
  $\bara _p$-module of rank $e$. Moreover $J \, Q_i \simeq Q_{i + 1}$. 
\end{prop}
\textbf{Proof.}  
  Since $x,y$ skew-commute, one readily computes that $Q_i$ is
  generated by $1, y, \ldots, y^{e - 1}$ as a $k
  [[z]]$-module. We wish to show they are free generators so suppose this is not the case and that $\text{rank}\ Q_i <e$. Now $\bara_p$ is a degree $e$ order so $Q_i$ must be torsion and $\xi_i$ must be invertible in $k(\bara_p)$. However, a direct computation shows $\xi_i^e = 0$, and this contradiction shows that $Q_i$ is torsion-free of rank $e$.
  
From lemma \ref{I-adicfiltrationonP} we know $y \, Q_i = J \, Q_i$ so is the cyclic module generated by $\bar{y}:=y + \bara_p \xi_i$. Now $\xi_{i+1} y = \z y \xi_i$ so $\xi_{i+1}$ annihilates $\bar{y}$ and there is a natural surjection $Q_{i+1} \lm JQ_i$. It is an isomorphism since $Q_{i + 1}$ and $J \, Q_i$ are free $k [[z]]$-modules of the same
  rank.
\vspace{2mm}

\begin{prop}
  Let $Q$ be a torsion-free $\bara _p$-module of rank $e$. Then $Q \simeq
  Q_i$ for some $i = 1, \ldots, e$.
\end{prop}
\textbf{Proof.}  
Note that left multiplication by $x$ or $y$ induces a $k [[z]]$-module
isomorphism $Q \to J \, Q$. Hence $y^{- 1} x$ induces a $k[[z]]$-module automorphism $\varphi$ of $Q$ and consequently, also of $Q/JQ$. Let $\bar{v} \in Q / J \,
Q$ be an eigenvector of $\varphi$. Its eigenvalue must have the form $\zeta^{i - (e - 1) / 2}$ since $(y^{-1}x)^e =  \zeta^{- e (e - 1) / 2}$. We lift $\bar{v}$ to an eigenvector $v \in Q$ so that $\xi_i v = 0$. By Nakayama's lemma and lemma~\ref{I-adicfiltrationonP}, $v$ generates $Q$ so there is a surjection $\psi:Q_i \lm Q$. Rank considerations show $\psi$ is an isomorphism.
\vspace{2mm}

We define the order $\bara _{1 p}$ as the subring of $\mathrm{End}_{k [[z]]}
(Q_e)$ consisting of all endomorphisms which stabilise $J^i Q_e$ for all $i =
1, \ldots, e - 1$. Since $J^i Q_e$ is an $\bara _p$-module for all $i
\geqslant 0$, we have $\bara _p \subseteq \bara _{1 p}$. The decomposition
\begin{eqnarray*}
  Q_e & = & k [[z]] \bar{x}^{e - 1} \oplus k [[z]] \bar{x}^{e - 2} \oplus
  \cdots \oplus k [[z]] \bar{x} \oplus k [[z]] \bar{1}
\end{eqnarray*}
induces an isomorphism $\mathrm{End}_{k [[z]]} (Q_e) \simeq k [[z]]^{e \times
e}$. By lemma \ref{I-adicfiltrationonP}, we have $J^i Q_e = x^i Q_e$, so we
can identify $\bara _{1 p}$ as the subring
\begin{eqnarray}
  \left(\begin{array}{cccc}
    k [[z]] & \cdots & \cdots & k [[z]]\\
    (z) & \ddots & \ddots & \vdots\\
    \vdots & \ddots & \ddots & \vdots\\
    (z) & \cdots & (z) & k [[z]]
  \end{array}\right) & \subseteq & k [[z]]^{e \times e} . 
  \label{standardheredorder}
\end{eqnarray}
It is clear from this description that $\bara _{1 p}$ is an hereditary order.
Furthermore, we can identify $\bara _p$ as the subring of $\bara _{1 p}$
generated by the following matrices
\begin{eqnarray}
  x = \left(\begin{array}{ccccc}
    0 & 1 & 0 & \cdots & 0\\
    \vdots & \ddots & \ddots & \ddots & \vdots\\
    \vdots & \ddots & \ddots & \ddots & 0\\
    0 & \ddots & \ddots & \ddots & 1\\
    z & 0 & \cdots & \cdots & 0
  \end{array}\right) & \mathrm{\mathrm{and}} & y = \left(\begin{array}{ccccc}
    0 & a_2 & 0 & \cdots & 0\\
    \vdots & \ddots & \ddots & \ddots & \vdots\\
    \vdots & \ddots & \ddots & \ddots & 0\\
    0 & \ddots & \ddots & \ddots & a_e\\
    a_1 z & 0 & \cdots & \cdots & 0
  \end{array}\right)  \label{embedding}
\end{eqnarray}
where $a_i = \zeta^{i + (e - 1) / 2}$. Hence the matrix $(\xi_{rs})$ corresponding to $\xi_j$ is also zero away from the $s-r \equiv 1 \mod e$ diagonal and $\xi_{e-j-1,e-j} = 0$. Let $\varepsilon_1, \ldots,
\varepsilon_e$ denote the standard matrix idempotents in $\bara _{1 p}$.

\begin{prop}
  \label{idempotentQ}We have an isomorphism of $\bara _{1 p}$-modules
  $\bara _{1 p} \varepsilon_j \simeq Q_{e - j}$.
\end{prop}
\textbf{Proof.}  
  Note that the composition $\varphi : \bara _p \hookrightarrow \bara _{1 p}
  \to \bara _{1 p} \varepsilon_j$ is surjective. Now $\xi_{e -
  j} \varepsilon_j = 0$, so $\bara _p \xi_{e - j} \subseteq \ker \,
  \varphi$ and we obtain an isomorphism $\bara _p / \bara  \xi_{e - j} = Q_{e
  - j} \to \bara _{1 p} \varepsilon_j$ of $\bara _p$-modules.
\vspace{2mm}


\begin{prop}
  \label{Ap1}Let $\bara _p \left\langle x^{- 1} y \right\rangle$ be the
  smallest subalgebra of $k ( \bara _p)$ containing $\bara _p$ and $x^{- 1}
  y$. Then we have $\bara _{1 p} = \bara _p \left\langle x^{- 1} y
  \right\rangle$. Moreover, as $k$-vector spaces, we have
  \begin{eqnarray}
    \bara _{1 p} & = & \bara _p \oplus \left( \bigoplus k \, x^i y^j z^{- 1}
    \right)  \label{kvsdecocmp}
  \end{eqnarray}
  where the direct sum ranges over $0 < i, j < e$ and $i + j \geqslant e$. In
  particular, $\dim_k ( \bara _{1 p} / \bara _p) = e (e - 1) / 2$.
\end{prop}
\textbf{Proof.}  
Firstly $x^{- 1} y$ is a diagonal matrix by (\ref{embedding}), so $\bara _p
\left\langle x^{- 1} y \right\rangle \subseteq \bara _{1 p}$. In fact, the entries of 
$\zeta^{- (e - 1) / 2}x^{- 1} y$ are the $e$ $e$-th roots of unity so $k[x^{- 1} y]$ is the subalgebra of $\bara_{1p}$ consisting of $k$'s along the diagonal in the matrix form (\ref{standardheredorder}). This shows that $\bara _p \left\langle x^{- 1} y \right\rangle$ contains all the standard matrix idempotents $\varepsilon_1, \ldots, \varepsilon_e$. Now
  $\varepsilon_{l - m} x^m \varepsilon_l$ generate $\varepsilon_{l - m}
  \bara _{1 p} \varepsilon_l$ as a $k [[z]]$-module, so $\bara _{1 p} =
  \bara _p \left\langle x^{- 1} y \right\rangle$.
  
  Note that $(x^{- 1} y)^i$ is a scalar multiple of $x^{e-i} y^i z^{- 1}$,
  so
  \begin{eqnarray*}
    \bara _{1 p} & = & \bara _p + 
\left( \bigoplus_{\begin{smallmatrix}i + j \geqslant e \\ 0 \leq i,j < e\end{smallmatrix}} 
k[[z]] \,x^i y^j z^{- 1} \right) .
  \end{eqnarray*}
But $\bara_p = \bigoplus_{0 \leq i,j < e} k[[z]] x^iy^j$ so the proposition follows. \ 
\vspace{2mm}

\begin{rem}
  One can show, using the criteria in [HN, Section 0.1], that the local ring
  $\bara _p$ does not have finite representation type if $e > 3$. This is in
  contrast with the situation in section \ref{sfibretan}. 
\end{rem}

Define the order $\bara _1$ on $F$ by $\bara _1 |_U = \bara  |_U$ where $U = F -
p$ and $( \bara _1)_p = \bara _{1 p}$. Since $\bara $ is maximal away from
$p$, it is clear that $\bara _1$ is an hereditary order totally ramified at $p$. This puts us in a good position to use the results of sections~\ref{sheredover} and \ref{sratgood}. We wish first to show that
$H^1 ( \bara _1) = 0$ by invoking lemma~\ref{lbdhered}. As in section~\ref{sfibretan}, we will need to carry out some local computations of ext groups and endomorphism rings in preparation for this.

Let $K_i = \bara _p \xi_i$,  $K_{1 i} = \bara _{1 p} \xi_i$ and 
$J_1$ be the Jacobson radical of $\bara _{1 p}$. Note that $J_1 = x
\bara _{1 p} = \bara _{1 p} x$.

\begin{prop}
  \label{propAp1}We have a natural isomorphism $\bara _{1 p} \simeq Q_1 \oplus
  \cdots \oplus Q_e$ and under this isomorphism, $K_{1 i}$ embeds as
  $\bigoplus_{j \neq i} J_1 Q_j \subseteq \bigoplus_{j = 1}^e Q_j$. Furthermore, $K_i = K_{1i} \cap \bara_p$.
\end{prop}
\textbf{Proof.}  
  The first isomorphism follows from $\bara _{1 p} = \bigoplus_{j = 1}^e
  \bara _{1 p} \varepsilon_j$ and proposition \ref{idempotentQ}. Note that
  $\varepsilon_j \xi_i = \xi_i \varepsilon_{j + 1}$ so $K_{1 i} =
  \bigoplus_{j = 1}^e \bara _{1 p} \xi_i \varepsilon_{j + 1}$. Moreover
  $\xi_i \varepsilon_{j + 1}$ is a scalar multiple of $x \varepsilon_{j + 1}$
  if $j \neq e - i - 1$, and $\xi_i \varepsilon_{e - i} = 0$. Since
  $\bara _{1 p} x = J_1$, we have $K_{1 i} = \bigoplus_{j \neq e - 1 -
  i} J_1 \bara _{1 p} \varepsilon_{j + 1}$. Finally, we see that $K_i = \ker(\bara_p \lm Q_i) \supseteq K_{1i} \cap \bara_p$ while the reverse inclusion is clear.
\vspace{2mm}

We omit the proofs of the next two elementary results.

\begin{lemma}
  \label{ext1SiPj}Let $S_i$ denote the simple $\bara _{1 p}$-module $Q_i /
  J_1 Q_i$. Then
  \begin{eqnarray*}
    \dim_k \mathrm{Ext}^1_{\bara _{1 p}} (S_i, Q_j) & = & \delta_{i + 1, j} .
  \end{eqnarray*}
\end{lemma}

\begin{lemma}
  \label{homlemma}Let $M$ be a torsion-free $\bara _p$-module and $N$ be a
  torsion-free $\bara _{1 p}$-module. Then there is a natural isomorphism
  \begin{eqnarray*}
    \mathrm{Hom}_{\bara _{1 p}} ( \bara _{1 p} M, N) & \simeq &
    \mathrm{Hom}_{\bara _p} (M, N) .
  \end{eqnarray*}
\end{lemma}

\begin{prop}
  \label{extcomputation}With the above notation, we have
  $\mathrm{Ext}^1_{\bara _p} (Q_j, Q_j) \simeq k$ and $\mathrm{Ext}^1_{\bara _p}
  (Q_j, K_j) \simeq k$ for any $j$.
\end{prop}
\textbf{Proof.}  
  We first show that $\mathrm{Ext}^1_{\bara _p} (Q_j, Q_j) \simeq k$. Firstly
  $\mathrm{Ext}^1_{\bara _p} (Q_j, Q_j) \simeq \mathrm{coker}
  (\mathrm{Hom}_{\bara _p} ( \bara _p, Q_j) \to
  \mathrm{Hom}_{\bara _p} (K_j, Q_j))$ and by lemma \ref{homlemma}, this is
  isomorphic to $\mathrm{coker} (\mathrm{Hom}_{\bara _{1 p}} ( \bara _{1 p},
  Q_j) \to \mathrm{Hom}_{\bara _{1 p}} (K_{1 j}, Q_j))$. This shows
  that
  \begin{eqnarray*}
    \mathrm{Ext}^1_{\bara _p} (Q_j, Q_j) & \simeq & \mathrm{Ext}^1_{\bara _{1
    p}} \left( \bara _{1 p} / K_{1 j}, Q_j \right) .
  \end{eqnarray*}
  By proposition \ref{propAp1}, we have $\bara _{1 p} / K_{1 j} \simeq Q_j
  \oplus \bigoplus_{i \neq j} S_i$ and the result follows from lemma
  \ref{ext1SiPj}.
  
  To compute $\mathrm{Ext}^1_{\bara _p} (Q_j, K_j)$, we use the following exact
  sequence
  \[ \mathrm{Hom}_{\bara _p} (Q_j, \bara _p) \overset{\rho_j}{\longrightarrow}
     k [[z]] \longrightarrow \mathrm{Ext}^1_{\bara _p} (Q_j, K_j)
     \longrightarrow \mathrm{Ext}^1_{\bara _p} (Q_j, \bara _p) . \]
The right most term is zero since $\bara _p$ has injective dimension one,
while $Q_j$ is torsion-free and hence is a first syzygy. The image of $\rho_j$ contains $\rho_j(\mathrm{Hom}_{\bara_p}(Q_j, z \bara _{1 p})) = \mathrm{Hom}_{\bara_{1p}}(Q_j, zQ_j)=(z)$. Also, the morphism $\rho_j$
  cannot be surjective for otherwise, any preimage of the
  identity map $\mathrm{id}_{Q_j}$ in $\mathrm{Hom}_{\bara _p} (Q_j, \bara _p)$
  would split the projection map $\bara _p \to Q_j$. This is
  impossible since $\bara _p$ is local. Hence 
  $\mathrm{Ext}^1_{\bara _p} (Q_j, K_j) \simeq k [[z]] / (z) = k$.
\vspace{2mm}

\begin{prop}
  \label{lengthcokernel} $\mathrm{End}_{\bara _p} (K_{1 i}) =
\mathrm{End}_{\bara _{1 p}} (K_{1 i})$ is naturally isomorphic to the non-unital  subring of $\bara_{1p}$ obtained by replacing the $(e-i)$-th row and column of the matrix form in (\ref{standardheredorder}) with 0s. In particular, it is a (unital) subring of $k[[z]]^{(e - 1) \times (e - 1)}$ with colength $(e - 1) (e - 2) / 2$.
\end{prop}
\textbf{Proof.}  
  This follows from proposition \ref{propAp1} and the natural isomorphism 
$$\Hom_{\bara_{1p}}(Q_i,Q_j) \xrightarrow{\sim} \Hom_{\bara_{1p}}(J_1Q_i,J_1Q_j):\phi \mapsto \phi|_{J_1Q_i}.$$

\begin{prop}
  \label{EndApKi}There is an exact sequence
  \[ 0 \longrightarrow \mathrm{End}_{\bara _p} (K_i) \longrightarrow
     \mathrm{End}_{\bara _{1 p}} (K_{1 i}) \longrightarrow K_{1 i} / K_i
     \longrightarrow 0. \]
  Moreover, $K_{1 i} / K_i$ is a $k$-vector space of
  dimension $(e - 1) (e - 2) / 2$. 
\end{prop}
\textbf{Proof.}  
We first establish the exact sequence above by considering the following exact sequence
  \[ 0 \longrightarrow \mathrm{End}_{\bara _p} (K_i) \longrightarrow
     \mathrm{Hom}_{\bara _p} (K_i, K_{1 i}) \overset{\alpha}{\longrightarrow}
     \mathrm{Hom}_{\bara _p} (K_i, K_{1 i} / K_i) . \]
  Note that by lemma \ref{homlemma}, the middle term is isomorphic to
  $\mathrm{End}_{\bara _{1 p}} (K_{1 i})$. Now the surjection $\bara _p \to \bara _p
  \xi_i = K_i$ induces an injection
  \[ \beta : \mathrm{Hom}_{\bara _p} (K_i, K_{1 i} / K_i) \longrightarrow
     \mathrm{Hom}_{\bara _p} ( \bara _p, K_{1 i} / K_i) \simeq K_{1 i} / K_i .
  \]
so it suffices to show surjectivity of $\alpha \circ \beta$.

We know from proposition~\ref{propAp1} that $K_{1i}/K_i \simeq (K_{1i} + \bara_p)/\bara_p$ so viewing $K_{1i}$ as a right $\End_{\bara_{1p}}(K_{1i})$-module, we are reduced to showing 
\begin{equation}
\xi_i \End_{\bara_{1p}}(K_{1i}) + \bara_{1p} = K_{1i} + \bara_{1p}
\label{eendK1}
\end{equation}
Matrix computations using propositions~\ref{propAp1} and \ref{lengthcokernel} show that 
$$ \xi_i \End_{\bara_{1p}}(K_{1i}) = \bigoplus_{j \neq e-i-1} \eps_j K_{1i} =: K'_{1i}$$
where $\eps_j$ are the standard diagonal matrix idempotents.
Furthermore, $x\xi_i, \ldots, x^{e-1}\xi_i \in \bara_p \cap K_{1i}$ generate  $K_{1i}/K'_{1i}$ as a $k[[z]]$-module so equation~(\ref{eendK1}) follows. This establishes the exact sequence of the proposition.

We now compute $\dim_k K_{1 i} / K_i$. From proposition~\ref{propAp1} we see there is an exact sequence of the form
$$  0 \lm K_{1i}/K_i \lm \bara_{1p}/J \lm \bara_{1p}/(K_{1i} + J) \lm 0 .$$
Now the natural map $J \lm JQ_i$ is surjective so proposition~\ref{propAp1} also shows that $K_{1i} + J = J\bara_{1p}$. Thus $\dim_k \bara_{1p}/(K_{1i} + J) = e$. On the other hand $\dim_k \bara_{1p}/J = e(e-1)/2 + 1$ by proposition~\ref{Ap1} so $\dim_k K_{1 i} / K_i = (e-1)(e-2)/2$ as was to be shown. 
%
%


\begin{thm}
  \begin{eqnarray}
    \bara _1 & \simeq & \left(\begin{array}{cccc}
      \mathcal{O} & \cdots & \cdots & \mathcal{O}\\
      \mathcal{O} (- 1) & \ddots & \ddots & \vdots\\
      \vdots & \ddots & \ddots & \vdots\\
      \mathcal{O} (- 1) & \cdots & \mathcal{O} (- 1) & \mathcal{O}
    \end{array}\right)  \label{hereditaryorderAbar1}
  \end{eqnarray}
\end{thm}
\textbf{Proof.}  
We consider an exact sequence of the form 
\[ 0 \longrightarrow K \longrightarrow \bara  \longrightarrow P
   \longrightarrow 0 \]
where $P$ is a torsion-free $\bara $-module of rank $e$. 
Note that proposition~\ref{idempotentQ} shows that $P_p \simeq Q_i$ for some $i$, while 
$K_p \simeq K_i$ by lemma~\ref{lsameker}. Hence lemma~\ref{lbdhered} reduces the proof to showing that $\nu (P) \geqslant 2 - e$ where
  \begin{eqnarray*}
    \nu (P) & = & \chi ( \send _{\bara } (P)) - \chi (
    \send _{\bara } (K)) - \dim_k \mathrm{Ext}^1_{\bara _p} (Q_i,
    Q_i) + \dim_k \mathrm{Ext}^1_{\bara _p} (Q_i, K_i) .
  \end{eqnarray*}
The dimensions of the ext terms above cancel by proposition \ref{extcomputation}. To calculate $\chi ( \send _{\bara } (K))$, we
  let $\varphi : \send _{\bara } (K) \to \mathcal{B}$
  be an embedding into a maximal order $\mathcal{B}$. Then $\chi ( \send _{\bara } (K)) = \chi ( \mathcal{B}) - \mathrm{length} (\mathrm{coker} \;
  \varphi)$. Since $\chi ( \mathcal{B}) = (e - 1)^2$ we conclude from
  proposition \ref{EndApKi} and proposition \ref{lengthcokernel} that
  \begin{eqnarray*}
    \chi ( \send _{\bara } (K)) & = & (e - 1)^2 - (e - 1) (e - 2)
    \; = \; e - 1.
  \end{eqnarray*}
  Finally we have $\send _{\bara } (P) = \mathcal{O}_F$, so $\nu
  (P) = 2 - e$ and the theorem is proved.
\vspace{2mm}

We next show that $\bara $ is rationally filtered. Let $P_i$ be the $i$-th
column of (\ref{hereditaryorderAbar1}), so that $\chi (P_i) = i$, and let
$\bara ^i = \bara  \cap P_{> e - i}$ as in section \ref{filtration_on_B}.
This gives a filtration
\[ 0 < \bara ^1 < \cdots < \bara ^{e - 1} < \bara ^e = \bara  . \]
\begin{thm}
\label{Aisrationallyfiltered} 
The filtration $\{ \bara ^j \}_{j = 1}^e$ is rational and $F$ is a good fibre. 
\end{thm}
\textbf{Proof.}  
By theorem \ref{ttannoH1}, it suffices to prove that $\chi ( \bara ^j) = j$. We calculate $\chi(\bara^j)$ using the exact sequence
  \[ 0 \longrightarrow \bara ^j \longrightarrow \bara  \overset{\psi}{\longrightarrow}
     P_{\leqslant e - j} \longrightarrow C_j \longrightarrow 0
  \]
  where $C_j = \mathrm{coker} \; \psi$. Note that $C_j$
  is supported at $p$ so we may compute $\chi(C^j)$ locally. Let $Q =
  (P_{\leqslant e - j})_p$. Since $Q$ is a summand of $\bara _{1 p}$, we
  have $Q \simeq Q_{i_1} \oplus \cdots \oplus Q_{i_{e - j}}$ where the indices
  are pairwise distinct. Let $\theta = \psi \otimes_F \hat{\calo}_{F,p}: \bara _p \to Q$, and note that $\mathrm{im} \; \theta$
  generates $Q$ as an $\bara _{1 p}$-module. Let $\theta_s : J^s / J^{s + 1}
  \to J^s Q / J^{s + 1} Q$ denote the map induced by $\theta$. We need two lemmas.
\vspace{2mm}

\begin{lemma}
  \label{Jadiccalculation}The map $\theta_s$ is injective if $0 \leqslant s <
  e - j$. Furthermore, $\theta_{e - j - 1}$ is an isomorphism. 
\end{lemma}
\textbf{Proof.}  
We first show that $\theta_s$ is injective for $0 \leq s < e-j$. Let $\theta(1) = \sum_{m=1}^{e-j} a_m$ where $a_m$ is the appropriate generator of $Q_{i_m}$. Now $\{x^s, x^{s-1}y, \ldots, y^s\}$ gives a $k$-basis for $J^s/J^{s+1}$ so it suffices to show that their images in $J^sQ/J^{s+1}Q$ remain linearly independent. Now lemma~\ref{I-adicfiltrationonP} shows that $\{x^sa_1,\ldots,x^s a_{e-j}\}$ gives a basis of $J^sQ/J^{s+1}Q$ so the 
$$\sum_{m=1}^{e-j} x^{s-l}y^l a_m + J^{s+1}Q = 
\sum_{m=1}^{e-j} \zeta^{-l(i_m + (e-1)/2)} x^s a_m + J^{s+1}Q$$
are linearly independent for $l = 0, 1, \ldots , s$. This proves injectivity of $\theta_s$. Furthermore, $\dim_k J^{e-j-1}/J^{e-j} = e-j = \dim_k J^{e-j-1}Q/J^{e-j}Q$ so $\theta_{e-j-1}$ is an isomorphism. 

\vspace{2mm}

The following is standard so we omit the proof. 
\begin{lemma}
  \label{injectivetest}Let $\varphi : M' \to M$ be a morphism of
  finite length modules over an Artin ring $R$. Denote by $J$ the Jacobson
  radical of $R$. Suppose that for all $k \geqslant 0$, the induced map
  $\varphi_k : J^k M' / J^{k + 1} M' \to J^k M / J^{k + 1} M$ is
  injective, then $\varphi$ is injective.
\end{lemma}

\noindent
\textbf{Rest of proof of theorem~\ref{Aisrationallyfiltered}.}  
  By lemma~\ref{Jadiccalculation}, $\theta_{e - j - 1}$ is an
  isomorphism so $\theta |_{J^{e - j - 1}}$ is
  surjective by Nakayama's lemma. Thus $\theta$ induces a map of left modules $\theta' : \bara _p / J^{e - j
  - 1} \to Q / J^{e - j - 1} Q$ with $\mathrm{coker} \; \theta' = C_j$. Furthermore, $\theta'$ is injective by lemmas~\ref{injectivetest} and \ref{Jadiccalculation}, so 
  \begin{eqnarray*}
    \mathrm{length} (C_j) & = & \mathrm{length} (Q / J^{e - j - 1} Q) -
    \mathrm{length} ( \bara _p / J^{e - j - 1})\\
    & = & (e - j) (e - j - 1) - \frac{(e - j) (e - j - 1)}{2}\\
    & = & \frac{(e - j) (e - j - 1)}{2}
  \end{eqnarray*}
  and
  \begin{eqnarray*}
    \chi ( \bara ^j) & = & \chi ( \bara ) - \chi (P_{\leqslant e - j}) +
    \chi (C_j)\\
    & = & e - \frac{(e - j) (e - j + 1)}{2} + \frac{(e - j) (e - j - 1)}{2}\\
    & = & j.
  \end{eqnarray*}
\vspace{2mm}

\section{Ramified fibres}  \label{sramfibres}  

Let $A$ be a degree $e$ ruled order. Suppose now that $F$ is a ramification curve. We will only consider the case where $F$ is totally ramified, that is, the ramification index of $A$ at $F$ is  $e$. Because of this, we will not need to consider the condition that $c_2$ is minimal. Note that if $e$ is prime then any ramification curve is totally ramified.

This time $\bara:= A \otimes_Z \calo_F$ is no longer an order so we analyse it as follows. Let $\eta$ be the generic point of $F$ and $\mmm_{\eta} \triangleleft \calo_{Z,\eta}$ be the maximal ideal of the local ring at $\eta$. The localisation $A_{\eta}$ of $A$ at $\eta$ is a maximal order in a division ring so we have complete information about it. In particular, its Jacobson radical $J_{\eta}$ is a principal ideal satisfying $J_{\eta}^e = \mmm_{\eta}A_{\eta}$ and it has residue ring $A_{\eta}/J_{\eta}$ the degree $e$ field extension of $k(F)$ which is ramified at the two points $D \cap F$ where $D$ is the bisection that $A$ is ramified on. We define $J := A \cap J_{\eta}$. We need

\begin{prop}  \label{pJi}  
For $i \in \mathbb{N}$ we have, 
i) $J^i$ is an invertible bimodule, 
ii) $J_{\eta}^i = J^i \otimes_Z \calo_{Z,\eta}$ and, 
iii) $J^i = A \cap J_{\eta}^i$.
\end{prop}
\textbf{Proof.} Note that $J$ is saturated in $A$ in the sense that $A/J$ embeds in $A/J \otimes_Z \calo_{Z,\eta}$. Hence $J$ is reflexive and thus an invertible bimodule. This proves i). Also, $J_{\eta} = J \otimes_Z \calo_{Z,\eta}$ so ii) follows. Finally, we see that $J^i$ is reflexive so also saturated. Thus iii) follows from ii).
\vspace{2mm}

\begin{prop}  \label{presring}  
The residue ring $A/J \simeq \calo_G$ where $G$ is the smooth rational curve which is the $e$-fold cover of $F$ ramified at the two points of $D \cap F$.
\end{prop}
\textbf{Proof.} We know that $A/J \hookrightarrow A_{\eta}/J_{\eta}$ and that $A/J \otimes_F k(F) \simeq A_{\eta}/J_{\eta}$. Hence $A/J \simeq \calo_G$ where $G$ is a projective model for the field extension of $k(F)$ totally ramified at $D\cap F$ with ramification index $e$. It suffices to show that $G$ is smooth. This can be done by local computations using theorem~\ref{tlocalterm}. 
\vspace{2mm}

We need to introduce some notation from [AV]. Given a sheaf $\mathcal{F}$ on a scheme $X$, and an automorphism $\s$ of $X$ we obtain an $\ox$-bimodule $\mathcal{F}_{\s}$ whose left module structure is given by $\mathcal{F}$ but the right module structure is $\s^*\mathcal{F}$. 
\begin{lemma}  \label{lJmodJ2}  
The $\calo_G$-bimodule $J/J^2$ is isomorphic to $(\calo_G)_{\sigma}$ where $\s \in \Aut G$ is the automorphism given by ramification data.
\end{lemma}
\textbf{Proof.} Ramification theory tells us that $J_{\eta}/J^2_{\eta}\simeq k(G)_{\s}$. Now $(J/J^2)^e \simeq J^e/J^{e+1} = A(-F)/J(-F) \simeq \calo_G$. Hence $\deg J/J^2 = 0$ as a sheaf on $G$. This proves the lemma.
\vspace{2mm}

Even though $\bara$ is not an order, $\bara\otimes_F k(F)$ still satisfies hypothesis~\ref{arankesimples} so we are in a good position to apply the results of section~\ref{sratgood}. 
\begin{thm}  \label{tramfibre} 
$\bara$ is rationally filtered so in particular, the fibre $F$ is good. Furthermore, if $c\in C$ is the point corresponding to the fibre $F$, then the map $C' \lm C$ is totally ramified at $c$.
\end{thm}
\textbf{Proof.} We consider the filtration of $\bara$
$$ 0 < J^{e-1}/A(-F) < \ldots < J/A(-F) < \bara .$$
The factors are all isomorphic to $\calo_G \simeq \calo_F \oplus \calo_F(-1)^{e-1}$ so the filtration is rational and $F$ is good by theorem~\ref{tratgood}. Also $J^2/A(-F)$ is rationally filtered so we may  lift a non-zero global section of $J/J^2$ to $\delta \in H^0(J/A(-F))$. One sees that $H^0(\bara) = k[\delta]/(\delta^e)$. Thus $C'\lm C$ is totally ramified at $c$. 

\section{The non-commutative Mori contraction}  \label{smori} 

In this section, we collate the results of the previous sections to prove theorem~\ref{tmain}. 

\begin{prop}  \label{plast}  
Let $A$ be a ruled order with minimal $c_2$ in its Morita equivalence class. Suppose further that every fibre is good. Then
\begin{enumerate}
\item $C':=\sspec_C f_*A$ is a smooth projective curve. 
\item $C'$ is a component of the Hilbert scheme $\Hilb A$ and $\Psi: A \otimes_Z f^*f_*A \lm A$ is the universal rational curve.
\item The non-commutative Mori contraction is given by the ``morphism of algebras'' $f_*A \lm A$.  
\end{enumerate}
\end{prop}
\textbf{Proof.} Parts~i) and ii) are essentially a restatement of proposition~\ref{pgood}. The morphism $\psi:f_*A \lm A$ determines a pull-back functor $\psi^*: f_*A-\Mod \lm A-\Mod$ as follows. Given a quasi-coherent $f_*A$-module $L$ then $f^*L$ is an $f^*f_*A$-module so we may consider the $A$-module $\psi^*(L):= A\otimes_{f^*f_*A} f^*L$. Since the universal rational curve is given by $\Psi$, we see that this coincides with the Fourier-Mukai transform of theorem~\ref{tCN}ii). Part~iii) follows immediately.
\vspace{2mm}

To see theorem~\ref{tmain}, it suffices to show that all fibres are good. This has been checked in theorems~\ref{theredfibres}, \ref{ttannoH1}, \ref{Aisrationallyfiltered} and \ref{tramfibre}.

\vspace{5mm}
\textbf{\large References}

\begin{itemize}
\item [{[A86]}] M. Artin, ``Maximal orders of global dimension and 
                   Krull dimension two'', Invent. Math., \textbf{84}, (1986), p.195-222
\item[{[A92]}] M. Artin, ``Geometry of quantum planes'', {\em Contemp. Math.}, \textbf{124}, (1992), p.1-15
\item[{[AdJ]}] M. Artin, A. J. de Jong, ``Stable Orders over Surfaces'', manuscript available at \\ www.math.lsa.umich.edu/courses/711/ordersms-num.pdf
\item[{[AG]}] M. Auslander, O. Goldman, ``Maximal orders'', {\em Trans. of the AMS} \textbf{97}, (1960), p.1-24
\item[{[AM]}] M. Artin, D. Mumford, ``Some elementary examples of unirational varieties which are not rational'', {\em Proc. LMS} \textbf{3}, (1972), p.75-95
\item[{[ATV90]}] M. Artin, J. Tate, M. Van den Bergh, ``Some algebras associated to automorphisms of curves'', Grothendieck Festschrift, Birkh\"auser, Basel (1990), p.33-85
\item[{[ATV91]}] M. Artin, J. Tate, M. Van den Bergh, ``Modules over regular algebras of dimension 3'', {\em Invent. Math.} \textbf{106}, (1991), p.335-89 
  \item  [{[AV]}] M. Artin, M. Van den Bergh,  ``Twisted Homogeneous Coordinate Rings'',
       {\em J. of Algebra}, \textbf{133}, (1990) , p.249-271
  \item  [{[AZ94]}] M. Artin, J. Zhang, ``Noncommutative projective schemes'', {\em Adv. Math.} \textbf{109} (1994), p.228-87
  \item  [{[AZ01]}] M. Artin, J. Zhang,  ``Abstract Hilbert Schemes'', {\em Alg. and Repr. Theory} \textbf{4} (2001), p.305-94
\item[{[C05]}] D. Chan, ``Splitting bundles over hereditary orders'', {\em Comm. Alg.} \textbf{333} (2005), p.2193-9
\item [{[Chan]}] K. Chan, ``Log terminal orders are numerically rational'', http://arxiv.org/abs/0911.5545
\item[{[CHI]}] D. Chan, P. Hacking, C. Ingalls, ``Canonical singularities for orders over surfaces'', {\em Proc. LMS} \textbf{98}, (2009), p.83-115
\item [{[CI]}] D. Chan, C. Ingalls, ``Minimal model program for orders over surfaces'' Invent. Math. \textbf{161} (2005) p.427-52
\item[{[CK03]}] D. Chan,  R. Kulkarni, ``Del Pezzo Orders on Projective Surfaces'', Advances in Mathematics  \textbf{173}  (2003), p.144-177
\item [{[CK05]}] D. Chan, R. Kulkarni, ``Numerically Calabi-Yau orders on surfaces''  {\em J. London Math. Soc.} \textbf{72}  (2005),  p.571-584
\item [{[CK11]}] D. Chan, R. Kulkarni, ``Moduli of bundles on exotic del Pezzo orders'', {\em Amer. J. Math.} \textbf{133}, (2011), p.273-93
\item [{[CN]}] D. Chan, A. Nyman, ``Non-commutative Mori contractions and $\PP^1$-bundles'', \\ http://arxiv.org/abs/0904.1717
\item[{[dJ]}] A.J. de Jong, ``The period-index problem for the Brauer group of an algebraic surface'',  {\em Duke Math. J.}  \textbf{123},  (2004),  p.71-94
\item[{[Ler]}] B. Lerner, PhD thesis UNSW, in preparation
\item[{[HN]}] H. Hijikata and K. Nishida, ``Primary orders of finite representation type'', {\em  J. Algebra}, \textbf{192} (1997) p.592-640
\item[{[Pot]}] J. le Potier, ``Lectures on vector bundles'', Cambridge studies in advanced mathematics \textbf{54}, Cambridge University Press, Cambridge, (1997)
\item[{[Ram]}] M. Ramras, ``Maximal orders over regular local rings of dimension two'' {\em Trans. of the AMS}, \textbf{142} (1969) p.457-79
\item[{[RV]}] I. Reiten, M. Van den Bergh, ``Two-dimensional Tame and 
          Maximal Orders of Finite Representation Type'', Mem. of the AMS 
          Vol. \textbf{80} No. 408, July 1989
\end{itemize}

\end{document}